\documentclass[a4paper,12pt]{article}

\usepackage{amsfonts,amssymb,amsbsy,amsmath,amsthm,mathrsfs,enumerate,verbatim}
\usepackage{pdfsync}

\usepackage{graphics} 
\usepackage{epsfig} 
\usepackage{graphicx}
\usepackage{algorithmic,algorithm}
\usepackage{epstopdf}
\usepackage[colorlinks=true]{hyperref}
\hypersetup{urlcolor=blue, citecolor=red}

\usepackage{amsfonts}
\usepackage{amsmath}
\usepackage{amssymb}
\usepackage{graphicx}
\usepackage{epsf}

\usepackage{mathptmx}       
\usepackage{helvet}         
\usepackage{courier}        
\usepackage{type1cm}        
\usepackage{makeidx}         
\usepackage{graphicx}        
\usepackage{multicol}        
\usepackage[bottom]{footmisc}
\usepackage{amsfonts}
\usepackage{amsmath}
\usepackage{amssymb}
\usepackage{amsfonts}
\usepackage{graphicx}
\usepackage{epsf}
\usepackage{epsfig}

\usepackage{color}
\usepackage{afterpage}
\usepackage{verbatim}
\usepackage{algorithmic,algorithm}

\makeindex             

\begin{document}

\title{Stability and convergence analysis of a domain decomposition FE/FD method for the Maxwell's equations in time domain}


\author{M. Asadzadeh$^*$ and L. Beilina
  \thanks{Department of Mathematical Sciences, Chalmers University of Technology and University of Gothenburg, SE-42196 Gothenburg, Sweden, e-mail: \texttt{\ mohammad@chalmers.se, larisa@chalmers.se}}
}

\date{}

\maketitle

\begin{abstract}

  Stability and convergence analysis for the domain decomposition finite element/finite difference (FE/FD) method developed in  \cite{BGrote, BMaxwell}   is presented. The analysis is designed for  semi-discrete  finite element scheme for the time-dependent Maxwell's equations. The explicit finite element schemes in different settings of the spatial domain are constructed and  domain decomposition algorithm  is formulated. Several numerical examples   validate   convergence rates obtained in the theoretical studies.
  
\end{abstract}

\section{Introduction}

New computational techniques meet the needs of
industry in developing efficient computational methods to simulate
partial differential equations (PDEs).  Especially,
for simulations in higher dimensions 
and large computational domains.  In this regard, certain
{\sl  domain decomposition } (DD)  method 
 leading to efficient schemes, in numerical investigations, 
 has gained a lot of interest in numerical analysis community. 
This variant of the DD method  is the
subject of our current, and some related and ongoing, research. 
This type of schemes was previously studied, e.g. in
\cite{Chan, Tos} for others than studied here problems.

 The present work is a further development of the
 DD hybrid finite element/finite difference (FE/FD)
 method for time-dependent Maxwell's equations for
 electric field in non-conductive media,  studied in \cite{BGrote,  BMaxwell}.  
 A stable, time-domain (TD), DDM scheme 
 for Maxwell's equations was
 proposed in \cite{RB1,RB2}, and further verified in \cite{GEMS1}.
 This method uses FDTD scheme on the, structured, FD part of the 
 mesh, and edge elements on unstructured part.  In applications, because of 
 edge elements implementation, this method remains computationally expensive.

 The domain decomposition FE/FD method for time-dependent
 Maxwell's equations for electric field, assuming  constant dielectric
 permittivity function in a finite difference domain, was 
 considered in \cite{BGrote}. This assumption 
 simplifies the numerical schemes in both FE
 and FD domains and significantly reduces computational efforts for
 implementation of the whole DD method.  Modified
 numerical scheme, energy estimate and numerical verifications of this
 method was presented in \cite{BMaxwell}. However, the fully stability and
 convergence analysis with numerical implementations, in $L_2$- and $H^1$-norms of the 
 developed FE and FD schemes, are  
 not presented in the above studies. We fill this gap  in the present work.
 
 More specificaly, we present stability analysis for explicit schemes
 for both FEM and FDM in the DD hybrid FE/FD method.
 The DDM is constructed such that FEM and FDM
 coincide on the common, structured, overlapping layer between the two
 subdomains.  The resulting domain decomposition approach at the
 overlapping layers can be viewed as a FE scheme
 which avoids instabilities at the interfaces.  Similar to the DD approach of \cite{BMaxwell, BeilinaHyb}, 
 we decompose
 the computational domain such that FEM and FDM are used in different subdomains: FDM 
  in simple geometry and FE  in the subdomain
 where more detailed information is needed about the structure of
 this subdomain. This also allows application of adaptive FEM in such 
 subdomain, see, e.g.  \cite{BMedical, BGrote, BTKM1, BTKM2, BondestaB, TBKF1, TBKF2}.

 Reliability and convergence of the domain decomposition method, 
 studied in this work, are  evident for solution of coefficient
 inverse problems (CIPs) in $\mathbb{R}^{3}$, see,
 e.g. \cite{BTKM1, BTKM2, BondestaB, TBKF1, TBKF2}. For the case
 of CIPs, the computational domain is splitted into subdomains such
 that a simple discretization scheme can be used in a large region and
 more refined discretization scheme is applied in smaller, however
 more critical, part of the domain.  In most algorithms for
 solution of electromagnetic CIPs, to determine the dielectric
 permittivity function inside a computational domain, a qualitative
 collection of experimental measurements is necessary on it's boundary
 or in it's neighborhood. In such cases it is convenient to condsider the numerical solution of
 time-dependent Maxwell's equations in
 different subdomains with constant dielectric permittivity function
 in some subdomain and non-constant in the other ones.  For the
 time-dependent Maxwell's equations, the DD scheme
 of \cite{BMaxwell}, which is analyzed in the present work, is used
 for solution of different CIPs to determine the dielectric
 permittivity function in non-conductive media using simulated and
 experimentally generated data, see \cite{BMedical, 
   BTKM1, BTKM2, BondestaB, TBKF1, TBKF2}.

An outline of this paper is as follows.  In Section 2 we introduce the
mathematical model.  In Section 3 we briefly present the domain
decomposition FE/FD method and communication scheme between two methods.  In
Section 4 we describe the domain decomposition FE/FD method for
solution of Maxwell's equations and set up the finite element and
finite difference schemes.  Section 5 is devoted to the stability analysis. 
In Section 6, we derive optimal a priori error estimates in
 finite element method for the 
semi-discrete (spatial discretizations) problems.  Finally, in Section 7 we present numerical
implementations that justify the theoretical investigations of the
paper. In what follows, $C$ will be a generic constant independent of all parameters, unless otherwise specifically specified, and not necessarliy the same at each occurance.


\section{The mathematical model}

\label{sec:model}

The Cauchy problem for the 
  electric field $E\left( x,t\right)=\left( E_{1},E_{2},E_{3}\right) \left( x,t\right)$,  
$x \in\mathbb{R}^{3}$,  $t \in [0, T]$, of the Maxwell's equations, under the 
assumptions  that the dimensionless relative magnetic permeability of 
the medium is $\mu_r \equiv 1$ and 
the  electric  volume
  charges are zero, is given by 
\begin{equation}\label{E_gauge}
\begin{split}
 \frac{1}{c^2} \varepsilon_r(x) \frac{\partial^2 E}{\partial t^2} +  \nabla \times \nabla \times E  &= - \mu_0 \sigma(x) 
 \frac{\partial E}{\partial t}, \\
  \nabla \cdot(\varepsilon E) &= 0,  \\
\end{split}
\end{equation}
where, $\varepsilon_r(x) = \varepsilon(x)/\varepsilon_0$ and
$\sigma(x)$ are the dimensionless relative dielectric permittivity and
electric conductivity functions, respectively. $\varepsilon_0$, and 
$\mu_0$ are the permittivity and permeability of the free space,
respectively, and $c=1/\sqrt{\varepsilon_0 \mu_0}$ is the speed of
light in free space. In this paper we consider the problem \eqref{E_gauge} in non-conductive media, i.e. 
$\sigma\equiv 0$, and hence study the initial value problem 

\begin{equation}\label{E_gauge1}
\begin{split}
 \varepsilon(x) \frac{\partial^2 E}{\partial t^2} +  \nabla \times \nabla \times E  &= 0, 
 \qquad x \in \mathbb{R}^{3},\, \, t\in (0,T].\\
  \nabla \cdot(\varepsilon E) &= 0,  \\
  E(x,0) = f_0(x), \qquad \frac{\partial E}{\partial t}(x,0) &= f_1(x) ,~~ x \in \mathbb{R}^{3},\, \, t\in (0,T].
\end{split}
\end{equation}


To solve the problem (\ref{E_gauge1}) numerically, we consider it in a bounded domain
$\Omega\subset \mathbb{R}^{n}, \, n=2,3$  (instead of whole $\mathbb{R}^{n}$), 
with boundary $\partial \Omega$, and employ a split scheme on 
$\Omega$: a hybrid, finite element/finite
difference scheme, {\sl kind of }
 domain decomposition, developed in \cite{BGrote, BMaxwell} and summarized in
Algorithm 1. More specifically, we divide the computational domain
$\Omega$ into two subregions, $\Omega_{\rm FEM}$ and $\Omega_{\rm
  FDM}$ such that $\Omega = \Omega_{\rm FEM} \cup \Omega_{\rm FDM}$
and $ \Omega_{\rm FEM}$ is a subset of the convex hul of 
$\Omega_{\rm FDM}$.  The function $\varepsilon(x)$ 
is assumed to be constant in $\Omega_{\rm FDM}$, and bounded and smooth 
 in $\Omega_{\rm FEM}$. 
The communication ß
between $\Omega_{\rm FEM}$ and $\Omega_{\rm FDM}$ is arranged using an  
overlapping mesh structure through a two-element thick layer around 
$\Omega_{\rm  FEM}$ as shown by blue and green common boundaries in Figure
\ref{fig:F1}.  The blue boundary is outer boundary of $\Omega_{\rm
  FEM}$ and inner boundary of $\Omega_{\rm FDM}$. Similarly, the green
boundary is the inner boundary of $\Omega_{\rm FEM}$ from which the
solution is copied to the green boundary of $\Omega_{\rm FDM}$.

The key idea with such a decomposition is to be able  apply different
numerical methods in different computational domains.  For the
numerical solution of \eqref{E_gauge1} in $\Omega_{\rm FDM}$ we use the
finite difference method on a structured mesh. In $\Omega_{\rm FEM}$, we
use finite elements on a sequence of unstructured meshes $K_h =
\{K\}$, with elements $K$ consisting of triangles in $\mathbb{R}^2$
and tetrahedra in $\mathbb{R}^3$, both satisfying minimal angle condition.
This approach combines the flexibility of the finite elements and the
efficiency of the finite differences in terms of speed and memory
usage and fits well for reconstruction algorithms presented below. 

We   assume that for some known constant $d>1$, the function
$\varepsilon\in C^{2}\left( \mathbb{R}^{3}\right)$ 
satisfies
\begin{equation} \label{2.3}
\begin{split}
  \varepsilon(x) \in \left[ 1,d\right],\quad
  & ~\text{ for }x\in  \Omega _{\rm FEM}, \\
  ~  ~ \varepsilon(x) =1, \qquad  & ~\text{ for } x\in  \Omega _{\rm FDM}. 
\end{split}
\end{equation}
Conditions \eqref{2.3} on $\varepsilon$  
and the relation 
\begin{equation}\label{divfree}
\nabla \times \nabla \times E = \nabla (\nabla
\cdot E) - \nabla \cdot ( \nabla E), 
\end{equation}
 together with divergence free field $E$, make equations in 
(\ref{E_gauge1}) independent of each others in $ \Omega _{\rm FEM}$ and $\Omega _{\rm FDM}$ so that, in  $\Omega _{\rm FDM}$,  
 we just need to solve the system of wave equations: 
\begin{equation} \label{6.11}
  \frac{\partial^2 E}{\partial t^2}-\Delta E= 0, ~~~(x,t) \in  \Omega_2 \times ( 0,T]. 
\end{equation}


\textbf{Remark}

It is well known that, for stable implementation of the finite element
solution of Maxwell's equations, divergence-free edge elements are the
most satisfactory ones from a theoretical point of view 
\cite{Monk0,  Nedelec}.  However, the edge elements are less attractive for
solution of time-dependent problems, since a linear system of equations
should be solved at each time iteration step.  In contrary, P1
elements can be efficiently used in a fully explicit finite element
scheme with lumped mass matrix \cite{Cohen, delta, joly}.  It is also
well known that numerical solution of Maxwell's equations using nodal
finite elements is often unstable and results spurious oscillatory
solutions \cite{MP, PL}.  There are a number of techniques 
to overcome such instabilities, see, e.g.  \cite{Jiang1, Jiang2, Jin, div_cor, PL}.



In \cite{BR1, BR2}, a finite element analysis shows stability
  and consistency of the stabilized  finite element method for the
  solution of (\ref{E_gauge}) with $\sigma(x)=0$. In the current study
  we show stability and convergence for the combined FEM/FDM scheme,
  under the condition \eqref{2.3} on $\varepsilon$, where the
  stabilized FEM is used for the numerical solution of
  (\ref{E_gauge1}) in $\Omega_{\mbox{FEM}}$ and usual FDM discretization of
  \eqref{6.11} is applied in $\Omega_{\mbox{ FDM}}$.

  \textbf{Remark}
    
Here, we consider the case when $ E(x,t) = 0$  for $x\in\partial \Omega$. Further, we assume that 
 $ E(x,0) \in [H^1(\Omega)]^3$, $\frac{\partial E}{\partial t} (x,0) \in {\textbf H}(\mbox{div},\Omega)$ and 
 $\nabla \cdot (\varepsilon  E(x,0)) = \nabla \cdot (\varepsilon \frac{\partial E}{\partial t} (x,0)) = 0$. 
 Recall that we assumed non-conductive media:  
 $\sigma \equiv 0$. 
In the presence of electric conductivity, additional $\sigma$-terms appear in the equations. 
 They lead to more involved estimates and heavier implementations which we plan 
 to perform in a forthcoming study.
 
Hence, in this note we study the following initial boundary value problem: 
  
\begin{equation}\label{eq1} 
\left\{
\begin{array}{ll}
    \varepsilon \partial_{tt} E + \nabla \times  \nabla \times E  = 0 & \mbox{ in } \Omega \times (0, T), \\
    E(\cdot,0) = f_0(\cdot), \mbox{ and } \partial_t E (\cdot,0) = f_1(\cdot) & \mbox{ in } \Omega, \\
    E = 0 & \mbox{ on } \partial \Omega \times (0,T), \\
    \nabla \cdot (\varepsilon E) = 0 & \mbox{ in } \Omega.
  \end{array}
  \right. 
\end{equation}


\section{The structure of domain decomposition}

\label{sec:hyb}

\begin{figure}[tbp]
\begin{center}
\begin{tabular}{ccc}
    {\includegraphics[trim = 6.0cm 0.0cm 0.0cm 0.0cm, scale=0.21 ,clip=]{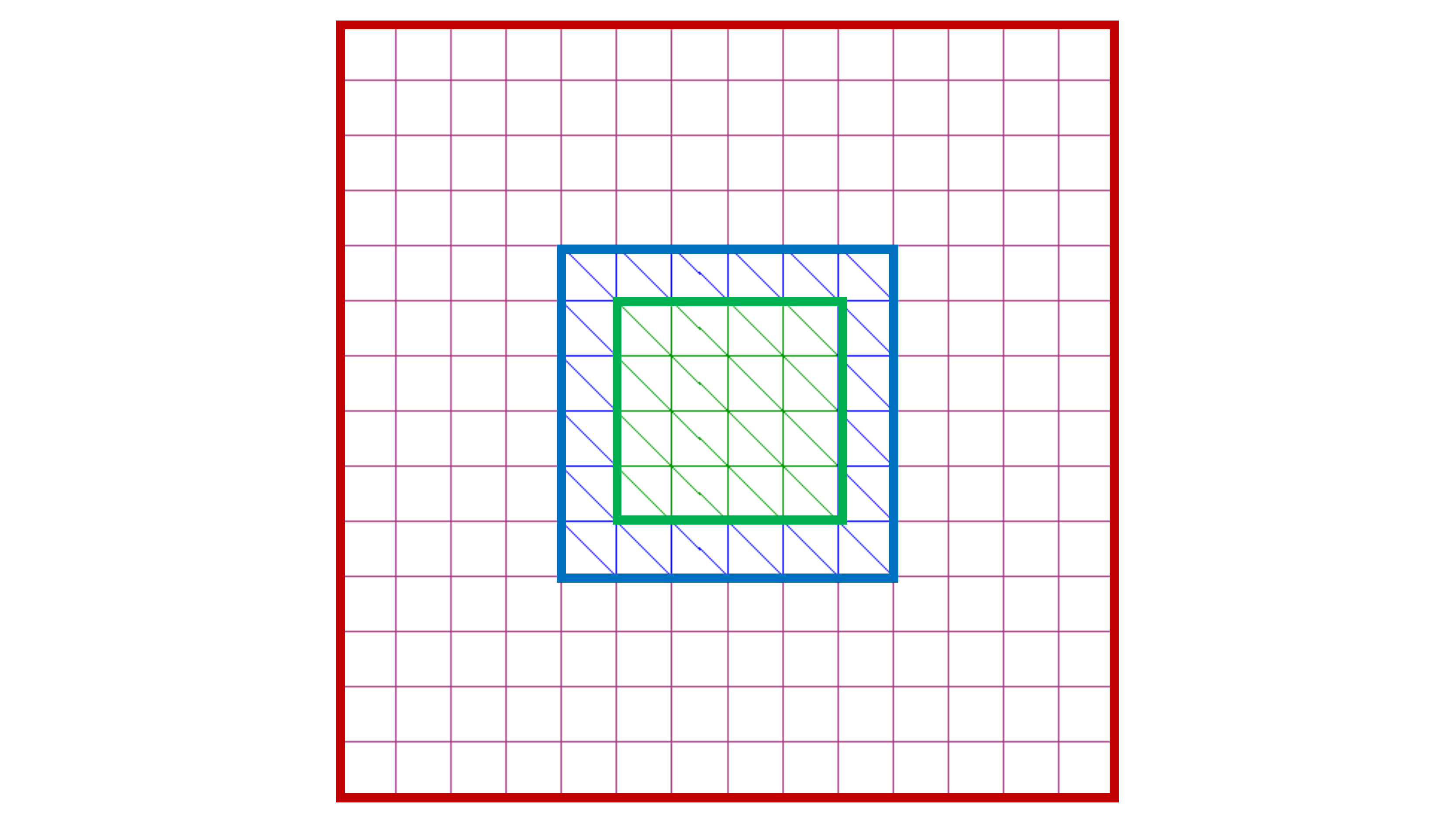}}
   &
   {\includegraphics[trim = 10.0cm 0.0cm 8.0cm 0.0cm, scale=0.21, clip=]{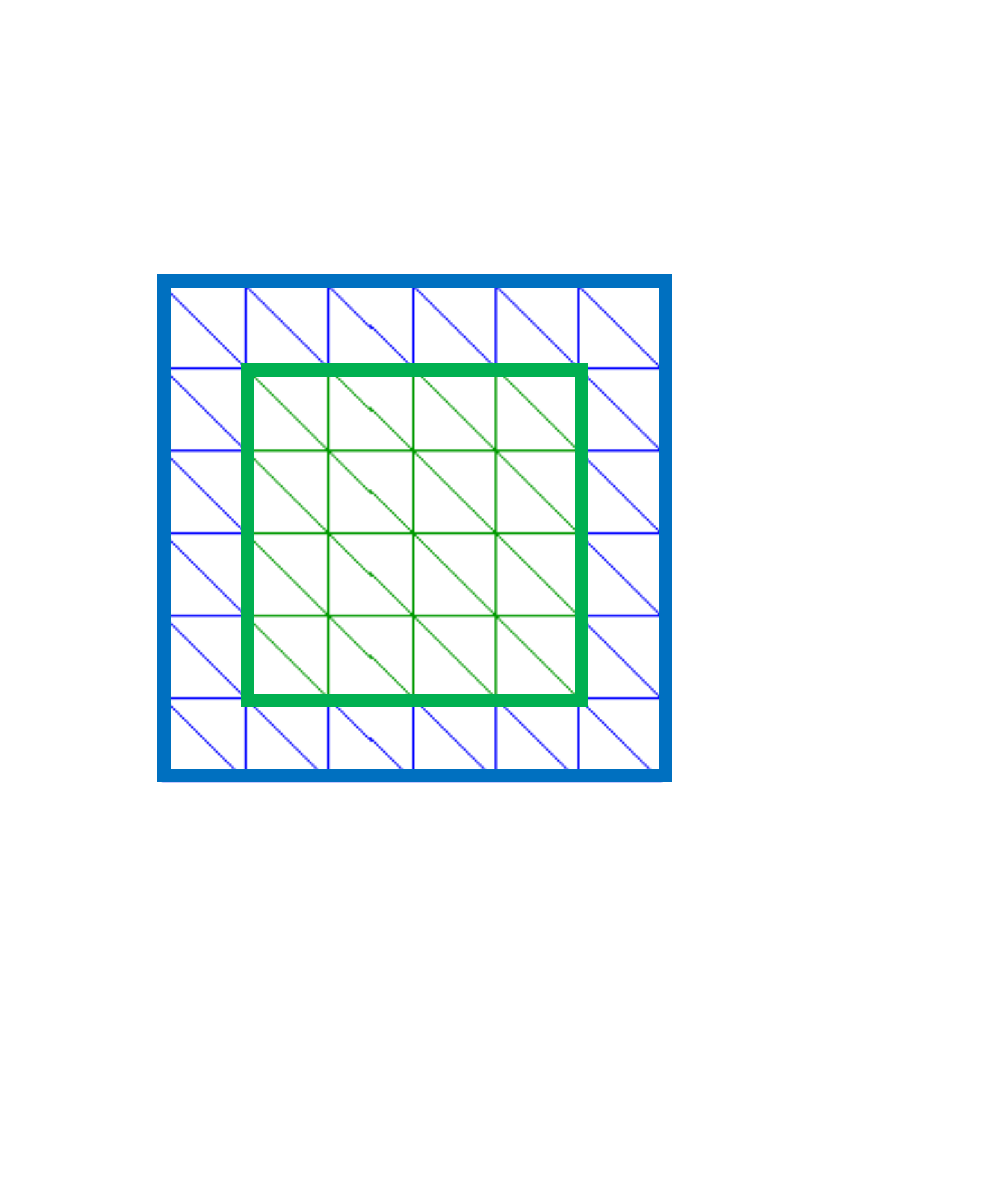}} 
 &
   {\includegraphics[trim = 6.0cm 0.0cm 0.0cm 0.0cm, scale=0.21,clip=]{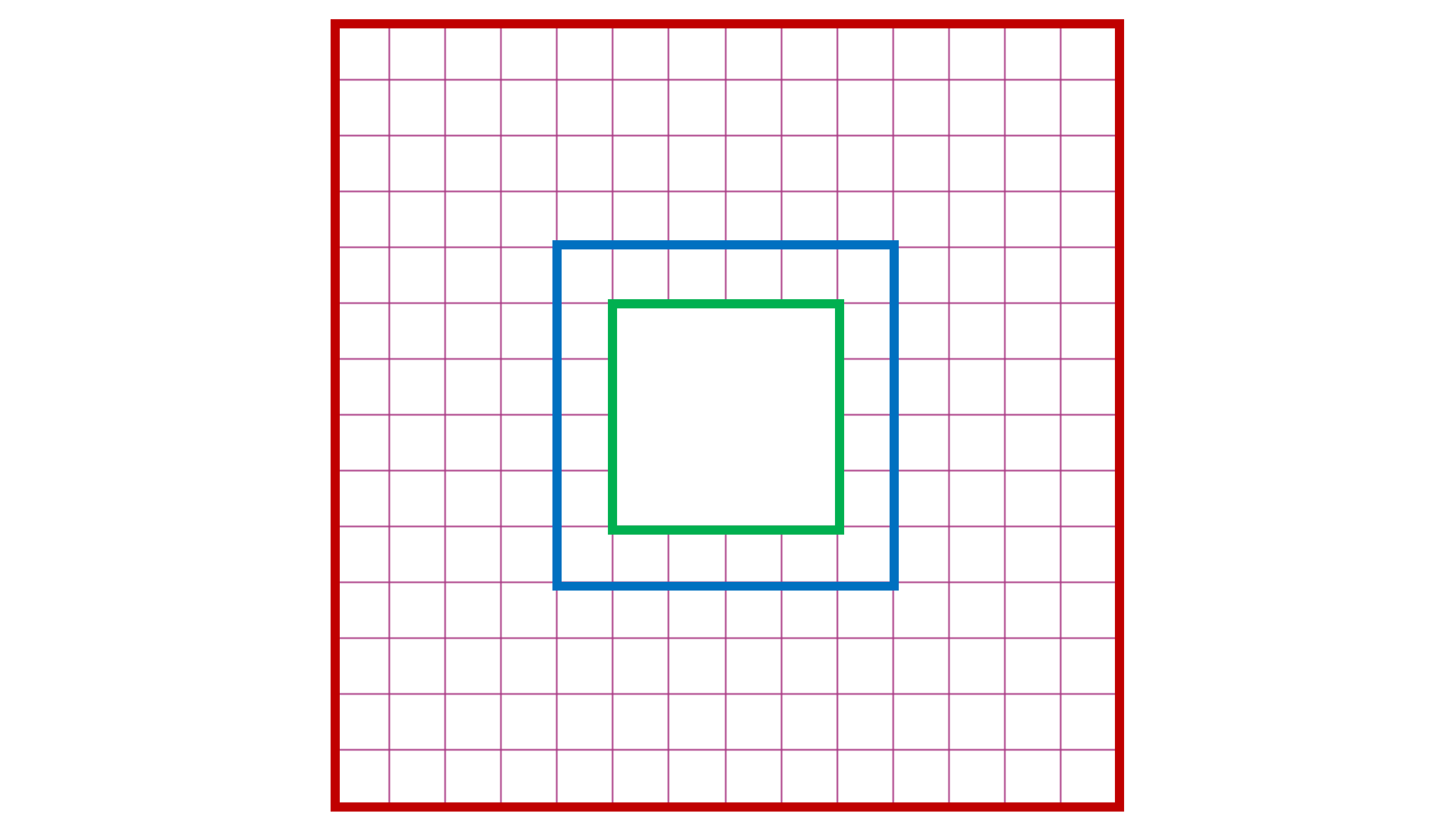}} \\  
 a) $ \Omega =  \Omega_{\rm{FEM} }\cup \Omega_{\rm{FDM}}$ &
 b) $\Omega_{\rm{ FEM}}$ & c) ~ $\Omega_{\rm{ FDM}}$
\end{tabular}
\end{center}
\caption{\small \emph{\, Domain decomposition and mesh discretization
    in $\Omega$.  The mesh of $\Omega$ is a combination of the
    quadrilateral finite difference mesh $\Omega_{\rm FDM}$ presented
    on c), and the finite element mesh $\Omega_{\rm FEM}$ presented on
    b). Domains $\Omega_{\rm FEM}$ and $\Omega_{\rm FDM}$ overlap by
    two layers of structured nodes such that they have common
    boundaries shown by green and blue colors.  }}
\label{fig:F1}
\end{figure}

We now describe the DD method between two domains, $\Omega_{FEM}$ and
$\Omega_{FDM}$, where the FEM is used for computation of the solution
in $\Omega_{FEM}$, and FDM is used in $\Omega_{FDM}$. Communication
between $\Omega_{FEM}$ and $\Omega_{FDM}$ is achieved letting
overlapping of both meshes across a two-element thick layer around
$\partial\Omega_{FEM}$ - see Figure \ref{fig:F1}.  The common nodes of both
$\Omega_{FEM}$ and $\Omega_{FDM}$ domains  belong to either of
the following boundaries (see Figure \ref{fig:F1}):
\begin{itemize}
\item Nodes on the blue boundary $\omega_{\rm o}$  - lie on the boundary $\partial \Omega_{FEM}$ of $\Omega_{FEM}$ and are interior to $\Omega_{FDM}$,
\item  Nodes on the green boundary  $\omega_{\diamond}$  -  lie on the inner boundary $\partial \Omega_{FDM}$ of
$\Omega_{FDM}$ and are 
 interior to $\Omega_{FEM}$.
\end{itemize}

\noindent 
Then the main loop in time for the explicit hybrid FEM/FDM scheme, that solves 
(\ref{E_gauge1}) 
 associates with appropriate boundary conditions, at each time step $k$  is described in 
 Algorithm 1 below:  

\vskip 0.3cm 

{
\begin{algorithm}[hbt!]
  \centering
  \caption{Domain decomposition process for   hybrid FE/FD scheme }
  \vskip 0.3cm 
  \begin{algorithmic}[1]
    \STATE   On  the mesh $\Omega_{FDM}$, where FDM is
  used, update the Finite Difference (FD) solution.

    \STATE  On the mesh $\Omega_{FEM}$, where FEM is used,  update the
    Finite Element (FE) solution.
    
    \STATE  Copy FE solution obtained at  nodes  $\omega_{\diamond}$ (nodes on the green boundary
   of Figure \ref{fig:F1})
    as a
  boundary condition on the inner boundary for the FD solution in $\Omega_{FDM}$.

  \STATE Copy FD solution obtained at nodes $\omega_{\rm o}$ (nodes  on the blue boundary
      of Figure \ref{fig:F1})
as a  boundary condition for the FE solution on $\partial \Omega_{FEM}$ of $\Omega_{FEM}$.

\STATE Apply boundary condition at $\partial \Omega$ at the red boundary of $\Omega_{FDM}$.

  \end{algorithmic}
\end{algorithm}
}

By  \eqref{2.3}, $\varepsilon =1$  at the overlapping nodes between $\Omega_{FEM}$ and 
$\Omega_{FDM}$. Thus,  FEM and FDM schemes coincide on the
common, structured, overlapping layer. Hence, we avoid
instabilities at interfaces.


\section{Derivation of computational schemes} \label{sec:domdec}

In this section we construct, combined, finite element-finite
difference schemes to solve the model problem \eqref{eq1}.  To do this, first  we present the finite 
element scheme  to solve
 \eqref{eq1} in entire $\Omega$. This induces 
 the finite element scheme in $\Omega_{FEM}$. Then, we derive the finite difference scheme in $\Omega_{FDM}$,   when  the
domain decomposition FE/FD structure is applied to solve  \eqref{eq1} in $\Omega$.


\textbf{Remark}

  The computational schemes derived in this section are explicit and therefore, for their  converegence, the CFL condition below (see ,e.g. \cite{BR1, BR2}) should be satisfied 
  \begin{equation}\label{CFL1}
\tau \leq \frac{h}{\eta}, ~~~\eta = C \sqrt{1 + 3 \| \varepsilon -1\|_\infty}. 
    \end{equation}
 Here $C $ is a mesh independent constant, 
 $\tau$ is the time step, and $h$ is the mesh size.


In the sequel, we denote the inner product of $[L^2(\Omega)]^M, M \in \{1,2,3\},$
by $(\cdot,\cdot)$, and the corresponding norm by $\parallel \cdot
\parallel$.
 The scalar  inner product  in $[L^2(\Omega_{FEM})]^M$  we denote by 
 $(\cdot,\cdot)_{\Omega_{FEM}}$, and the associated norm by $\| \cdot \|_{\Omega_{FEM}}$. 
Further,  we let $\partial \Omega_{FEM}$ to be the boundary of $\Omega_{FEM}$, $\partial \Omega_{FDM}$ the inner boundary of $\Omega_{FDM}$ and $\partial \Omega$  the outer boundary of $\Omega_{FDM}$.

\subsection{Finite element discretization in $\Omega$}
\label{sec:fem}

First we derive finite element scheme to solve the model  problem  \eqref{eq1}  in whole 
$\Omega$. 
Next,  we discretize $\Omega_{{FEM}_T} = \Omega_{FEM} \times (0,T)$ in two steps: 
(i) the spatial discretization using 
a partition $K_h = \{K\}$  of $\Omega_{FEM}$ into elements $K$, where 
$h=h(x)$ is a mesh function
 defined as $h |_K = h_K$, representing the local diameter of elements.
We also denote by  $\partial K_h = \{\partial K\}$ a partition of
 the boundary $\partial \Omega_{FEM}$ into boundaries  $\partial K$ of the elements $K$ such that, at least one of the vertices of these elements belong to  $\partial \Omega_{FEM}$.
(ii) As for temporal discretization, 
we let $J_{\tau}$ be a uniform partition of the time interval $(0,T)$ into $N$ 
 subintervals $J=(t_{k-1},t_k]$ of length $\tau=T/N.$
 As usual,  
 we also 
 assume a minimal angle condition on elements $K$ in $K_h$.
To formulate the  finite element method for \eqref{eq1} in $\Omega$,
   we introduce the finite element space $W_h^E(\Omega)$ for each component of the electric  field $E$ defined by
\begin{equation}
W_h^E(\Omega) := \{ w \in H^1(\Omega): w|_{K} \in  P_1(K),  \forall K \in K_h\}, \nonumber
\end{equation}
where $P_1(K)$ denote the set of piecewise-linear functions on $K$.
 Setting ${\textbf W_h^E(\Omega)} := [W_h^E(\Omega)]^3$ we define ${f_0}_h$ (resp. ${f_1}_h$) to be the 
 usual ${\textbf  W_h^E}$-interpolant of 
$f_0$ (resp. $f_1$) in \eqref{eq1}. 
   Also, because of the Dirichlet boundary data in 
\eqref{eq1} we need to choose the test function space as 
\begin{equation}\label{TestSpace}
{\textbf  W}_{h,0}^E:=\{{\mathbf v}\in {\textbf  W_h^E} | {\mathbf v}=0,\quad  \mbox{on }\quad \partial \Omega\}.
\end{equation}
Then, recalling  \eqref{divfree}, and the fact that 
$\nabla \cdot (\varepsilon E)=0,$ 
 the spatial semi-discrete problem 
 in $\Omega$ 
  reads: 

\emph{Find } $E_{h}\in {\mathbf W_h^E(\Omega)}$ \emph{ such  that } 
$\forall {\mathbf v} \in {\mathbf  W}_{h,0}^E(\Omega)$, 
\begin{equation}\label{eq6}
  \begin{array}{l}
    \left ( \varepsilon \partial_{tt} E_{h}, {\textbf v} \right ) + (\nabla E_h,\nabla {\textbf v})+ 
    (\nabla \cdot (\varepsilon E_h), \nabla \cdot {\textbf v} ) 
    -(\nabla \cdot E_h, \nabla \cdot {\textbf v}) \\
\qquad - \langle\partial_n E_{h}, {\textbf v}\rangle_{\partial \Omega} := 
\sum_{j=1}^5 T_j=0,  \\
\qquad   E_h(\cdot,0) = {f_0}_h(\cdot) \quad \mbox{ and } \,\, \partial_t E_h(\cdot,0) = {f_1}_h(\cdot) \quad 
   \mbox{ in } \Omega.  
  \end{array}
\end{equation}
We note that 
$$
T_3= (\nabla \cdot (\varepsilon E_h), \nabla \cdot {\textbf v} ) =
(\nabla\varepsilon \cdot E_h, \nabla \cdot {\textbf v} ) +
(\varepsilon \nabla\cdot E_h, \nabla \cdot {\textbf v} ), 
$$
implies that
\begin{equation}\label{eq6A1}
\begin{split}
T_3+T_4 &= (\nabla \cdot (\varepsilon E_h), \nabla \cdot {\textbf v} )-
(\nabla\cdot E_h, \nabla \cdot {\textbf v} ) \\
&= (\nabla\varepsilon \cdot E_h, \nabla \cdot {\textbf v} ) +
\Big( (\varepsilon -1)\nabla\cdot E_h, \nabla \cdot {\textbf v} \Big).
\end{split}
\end{equation}
Recalling \eqref{TestSpace}, vanishing test  functions at the boundary yields 
$T_5\equiv 0$, and 
 the final weak formulation for the semi-discrete problem in  $\Omega$ is: 
\emph{Find } $E_{h}\in {\textbf  W_h^E(\Omega)}$
\emph{ such  that} $\forall {\mathbf v} \in {\mathbf  W}_{h,0}^E(\Omega)$, 
\begin{equation}\label{eq6SEM}
   B_\Omega (E_h,  {\textbf v}):= 
   \left ( \varepsilon \partial_{tt} E_{h}, {\textbf v} \right ) +
    (\nabla E_h,\nabla {\textbf v})+   (\nabla \varepsilon \cdot E_h, \nabla \cdot {\textbf v} ) 
    +((\varepsilon-1)\nabla \cdot E_h, \nabla \cdot {\textbf v}) =0. 
\end{equation}
For the, reflexive, inhomogeneous boundary condition, see the FE scheme for $\Omega_{FEM}$.

To get fully discrete scheme for
 \eqref{eq1} we apply time discretization to \eqref{eq6} approximating $E_h(k\tau)$, denoted 
by $E_h^k$, where we use the central difference scheme, for $k=1,\ldots,N-1$:
\begin{equation}\label{eq7}
  \begin{array}{l}
    \left (\varepsilon \frac{E_h^{k+1} - 2 E_h^k + E_h^{k-1}}{\tau^2}, {\textbf v} \right) + (\nabla E_h^k, \nabla {\textbf v}) 
     + (\nabla \cdot (\varepsilon E_h^k), \nabla \cdot {\textbf v}) \\
      -(\nabla \cdot E_h^k, \nabla \cdot {\textbf v}) 
	= 0 \, \quad 
		 \forall {\textbf v} \in {\textbf  W_{h_0}^E(\Omega)},\\
		 {}
		\\
 \quad     E_h^0 = {f_0}_h\quad  \mbox{ and } \,\, E_h^1 = E_h^0 + \tau {f_1}_h \quad \mbox{ in } \Omega.
  \end{array}
\end{equation}  

Multiplying both sides of \eqref{eq7} by $\tau^2/\varepsilon$ we get,  
$ \forall {\mathbf v} \in {\textbf W}_{h,0}^E(\Omega),$ that 

\begin{equation}\label{eq8}
  \begin{array}{l}
    \left (E_h^{k+1} - 2 E_h^k + E_h^{k-1}, {\textbf v} \right)
    +\tau^2 (1/ \varepsilon \nabla E_h^k, \nabla {\textbf v})
 +\tau^2 (1/ \varepsilon \nabla \cdot (\varepsilon E_h^k), \nabla \cdot {\textbf v}) \\
    \qquad -\tau^2 (1/ \varepsilon \nabla \cdot E_h^k, \nabla \cdot {\textbf v}) =0, \\
    {}
    \\
     \qquad  E_h^0 = {f_0}_h \quad \mbox{ and } \quad E_h^1 =
          E_h^0 + \tau {f_1}_h,\qquad  \mbox{ in } \Omega.
  \end{array}
\end{equation}

Rearranging terms in \eqref{eq8} we get the following scheme:\emph{ Given the approximate initial data, ${f_0}_h$ and ${f_1}_h$}, 
\emph{find $E_{h}\in {\textbf  W_h^E(\Omega)}$  such  that $\forall {\mathbf v} \in {\mathbf  W_{h,0}^E(\Omega)}$ }, 
\begin{equation}\label{eq9}
  \begin{split}
   \left(E_h^{k+1}, {\textbf v} \right) & =   
 \left(2 E_h^k, {\textbf v} \right) - \left(E_h^{k-1}, {\textbf v} \right) 
     - \tau^2 (1/ \varepsilon \nabla E_h^k, \nabla {\textbf v})
  - \tau^2 (1/ \varepsilon \nabla \cdot (\varepsilon E_h^k), \nabla \cdot {\textbf v}) \\
  &+   \tau^2 (1/ \varepsilon \nabla \cdot E_h^k, \nabla \cdot {\textbf v})\\
  {}
  \\
 \,\,\,   \qquad  & E_h^0 = {f_0}_h \quad \mbox{ and } \,\,  E_h^1 =
          E_h^0 + \tau {f_1}_h \quad \mbox{ in } \Omega.
  \end{split}
\end{equation}

\subsection{Finite element discretization in $\Omega_{FEM}$}
\label{sec:omegafem}

To solve the model problem \eqref{eq1} via the domain decomposition
FE/FD method, we use the split $\Omega =\Omega_{FEM} \cup \Omega_{FDM}$, see Figure \ref{fig:F1}. 
Thus in $\Omega_{FEM}$ we use FEM to solve the equation
\begin{equation}\label{eq1fem}
  \begin{array}{ll}
    \varepsilon \partial_{tt} E + \nabla \times  \nabla \times E  = 0& \mbox{ in } \Omega_{FEM} \times (0, T), \\
    E(\cdot,0) = f_0(\cdot), \mbox{ and } \partial_t E (\cdot,0) = f_1(\cdot) & \mbox{ in } \Omega_{FEM}, \\
    \partial_n E =-\partial_t E = g & \mbox{ on } \partial \Omega_{FEM} \times (0,T), \\
    \nabla \cdot (\varepsilon E) = 0 & \mbox{ in } \Omega_{FEM}.
  \end{array}
\end{equation}
Here, $g$ is the restriction  of the solution obtained by the FDM in 
$\Omega_{FDM}$ to $\partial\Omega_{FEM}$ and therefore the test functions are 
not vanishing at the boundary and hence the term corresponding to $T_5$ in \eqref{eq6} will be appearing in 
the weak formulation. 

To formulate the  finite element method for \eqref{eq1fem} in $\Omega_{FEM}$, mimiking 
\eqref{eq6}, 
   we introduce the finite element space $W_h^E(\Omega_{FEM})$ for each component of the electric  field $E$ defined by
\begin{equation}
W_h^E(\Omega_{FEM}) := \{ w \in H^1(\Omega_{FEM}): w|_{K} \in  P_1(K),\quad  \forall K \in K_h\}. \nonumber
\end{equation}
 Setting ${\textbf W_h^E(\Omega_{FEM})} := [W_h^E(\Omega_{FEM})]^3$ we define 
 ${f_0}_h$, ${f_1}_h$, and $g_h$ to be the usual ${\textbf  W_h^E}$-interpolants of 
 $f_0$, $f_1$, and $g$, respectively, in $\Omega_{FEM}$. Then,  similar to the FE scheme for $\Omega$,   
 we get the following finite element scheme for $\Omega_{FEM}$:   
  {\sl Given ${f_0}_h$, ${f_1}_h$, and $g_h$, find 
$E_h\in  {\textbf W_h^E(\Omega_{FEM})}$ such that}
\begin{equation}\label{eq6B}
  \begin{split}
    \left ( \varepsilon \partial_{tt} E_{h}, {\textbf v} \right ) &+ (\nabla E_h,\nabla {\textbf v})+ 
    (\nabla \cdot (\varepsilon E_h), \nabla \cdot {\textbf v} ) 
    -(\nabla \cdot E_h, \nabla \cdot {\textbf v}) \\
 &=\langle g_{h}, {\textbf v}\rangle_{\partial \Omega_{FEM}} ,  \quad \forall  {\textbf v} \in  {\textbf W_h^E(\Omega_{FEM})} \\
 &{}
 \\
\qquad  & E_h(\cdot,0) = {f_0}_h(\cdot) \quad \mbox{ and } \,\, \partial_t E_h(\cdot,0) = {f_1}_h(\cdot) \quad 
   \mbox{ in } \Omega_{FEM}.  
  \end{split}
\end{equation}
 
 A corresponding fully discrete problem in $\Omega_{{FEM}_T}$ reads as follows: 
 
{ Given  ${f_0}_h$, ${f_1}_h$,  $g_h$, $E_h^k$, and $E_h^{k-1}$; find $E_h^{k+1}$ such that} 

\begin{equation}\label{eq9fem}
  \begin{split} 
    \left ( E_h^{k+1}, {\textbf v} \right)
  =  &  \left(2 E_h^k, {\textbf v} \right) - \left(E_h^{k-1}, {\textbf v} \right)
   - \tau^2 (1/ \varepsilon \nabla E_h^k, \nabla {\textbf v})\\
 & - \tau^2 (1/ \varepsilon \nabla \cdot (\varepsilon E_h^k), \nabla \cdot {\textbf v}) 
 +   \tau^2 (1/ \varepsilon \nabla \cdot E_h^k, \nabla \cdot {\textbf v}) \\
 &+\tau^2 \langle {g_h}/{\varepsilon}, {\textbf v} \rangle_{\partial \Omega_{FEM}}, 
 \qquad \forall {\textbf v} \in {\textbf W_h^E(\Omega_{FEM_T})}
\\
& 
\\
 \qquad E_h^0 = {f_0}_h & \quad E_h^1 =
          E_h^0 + \tau {f_1}_h\quad  \mbox{ in } \Omega_{{FEM}_T}.
  \end{split}
\end{equation}

\textbf{Remark}

Note that, in  \eqref{eq1fem}, Dirichlet boundary condition $E=g$  can be
  considered as well.

\subsection{Fully discrete FE scheme for the electric field  in $\Omega_T$}

\label{sec:discrete}

 We expand the functions $E_h$  in
terms of the standard continuous piecewise linear functions
$\{\varphi_i(x)\}_{i=1}^M$ in space as
\begin{equation}\label{eq10}
\begin{split}
E_h = \sum_{i=1}^M E_{h_{i}}(t) \varphi_i(x),
\end{split}
\end{equation}
 where $E_{h_{i}}(t)$ denote unknown coefficients at $t\in (0,T]$ and the spatial mesh point 
 $x_i \in K_h$. Then,  substituting $E_h$ of \eqref{eq10} in 
(\ref{eq9}), and setting 
${\mathbf v} = \sum_{j=1}^M \varphi_j(x)$,
we obtain the linear system  of equations:
\begin{equation} \label{femod1}
\begin{split}
 M  E^{k+1}_h &=   
2 M E^{k}_h  - M E^{k-1}_h  - \tau^2 G_1 E^{k}_h  -  \tau^2 G_2 E^{k}_h  \\
&+  
\tau^2 G_3 E^{k}_h + \tau^2  M_{\partial K} E^{k}_h. 
\end{split}
\end{equation}
Note that, unlike \eqref{eq9},  now the contributions at boundary of the element appear in 
$M_{\partial K}$. 
  Here, $M, M_1, M_2,  M_{\partial K}$ are the block mass matrices in space, $G_1,
  G_2, G_3$ are the block stiffness matrices in space, $E^k_h$ denote the nodal values
  of $E_h(\cdot,t_k)$, and $\tau$ is a uniform time step.
Now we define the mapping $F_K$ from the reference element $\hat{K}$ onto $K$ 
such that $F_K(\hat{K})=K$ and let $\hat{\varphi}$ be the piecewise
linear local basis function on $\hat{K}$ such
that $\varphi \circ F_K = \hat{\varphi}$.  Then, the explicit formulas
for the entries in system of equations (\ref{femod1}), for each element $K$, can be
written as:
\begin{equation}\label{eq11} 
\begin{split}
 M_{{i, j}}^{K}  &=    (\varphi_i(x) \circ F_K, \varphi_j(x)\circ F_K)_K, \\
 M_{{i, j}}^{\partial K}  &=  
 \langle\frac{1}{ \varepsilon}\partial_n \varphi_i(x) \circ F_K, \varphi_j(x) \circ F_K\rangle_{\partial K},\\
 {G_1}_{{i, j}}^{ K} & =(\frac{1}{\varepsilon} \nabla  \varphi_i \circ F_K, \nabla  \varphi_j \circ F_K)_K, \\
{G_2}_{{i, j}}^{ K} & =(\frac{1}{\varepsilon} \nabla \cdot (\varepsilon \varphi_i) \circ F_K, \nabla \cdot \varphi_j \circ F_K)_K, \\
{G_3}_{{i,j}}^{ K} & =(\frac{1}{\varepsilon} \nabla \cdot \varphi_i \circ F_K, \nabla \cdot \varphi_j \circ F_K)_K, 
\end{split} 
 \end{equation}
where $(\cdot,\cdot)_K$, and $\langle \cdot,\cdot\rangle_{\partial K}$, denote 
the $L_2(K)$, and  $L_2(\partial K)$, scalar products on $K$  and $\partial K$, respectively. Note that here 
 $\partial K$ is only the part of the boundary of element $K$ that lies at $\partial \Omega_{FEM}$.

To obtain fully explicit scheme we approximate $M$ with the lumped mass
matrix $M^{L}$, (see \cite{delta, joly, BR1} for the details corresponding to 
the Maxwell's system \eqref{E_gauge1}). 
Next, we multiply (\ref{eq11}) by $(M^{L})^{-1}$, and
get the following explicit, fully discrete method in $\Omega$: 
\begin{equation} \label{femod1n}
\begin{split}
  E^{k+1}_h &=   
2  E^{k}_h  -  E^{k-1}_h  - \tau^2 (M^{L})^{-1}  G_1 E^{k}_h  -  \tau^2 (M^{L})^{-1} G_2 E^{k}_h  \\
&+  
\tau^2  (M^{L})^{-1} G_3 E^{k}_h + \tau^2  (M^{L})^{-1}   M_{\partial K} E^{k}_h. 
\end{split}
\end{equation}

\subsection{Fully discrete scheme for the electric field  in $\Omega_{FEM_T}$}

\label{sec:discretefem}

 As in the fully discrete FE scheme
  \eqref{femod1} in $\Omega_T$,
we obtain  fully discrete FE scheme  in $\Omega_{{FEM}_T}$
in the domain decomposition setting: 
Expanding the  $E_h$ functions of $\Omega_{{FEM}}$ via the
  continuous piecewise linear functions in space as in \eqref{eq10}, 
  and then substituting them in \eqref{eq9fem}, (with 
  ${\mathbf v}= \sum_{j=1}^M\varphi_j(x), \, x\in  \Omega_{{FEM}}$, and 
  $x_i\in K_h\subset \Omega_{{FEM}}$ ), we  get the linear system of  equations:
\begin{equation} \label{femod2}
\begin{split}
 M E^{k+1}_h =   & 
2 M E^{k}_h  - M E^{k-1}_h  - \tau^2 G_1 E^{k}_h  -  \tau^2 G_2 E^{k}_h  \\
&+  
\tau^2 G_3 E^{k}_h + \tau^2  S_{\partial K} . 
\end{split}
\end{equation}
  Here, $M$ is the block mass matrices in space, restricted to $\Omega_{{FEM}}$, otherwise the
  same as in \eqref{femod1}, $G_1, G_2, G_3$ are the block stiffness matrices in
  space as in \eqref{femod1}, $S_{\partial K}$ is the assembled load
  vector , $E^k_h$ denote the nodal values of $E_h(\cdot,t_k)$, $\tau$
  is the time step.  All  quantities are for $\Omega_{{FEM}_T}$. 
  Defining the mapping $F_K$ for the reference
  element $\hat{K}$ in the mesh $K_h$ generated in $\Omega_{FEM}$ as
  in the previous section, the formulas for entries of all matrices in the system 
 \eqref{femod2} are the same as those  in \eqref{eq11}, and the entries of load
  vector are computed as
\begin{equation}
  S_{j}^{K}=\Big(  \frac{g_h}{ \varepsilon}, \varphi_j \circ F_K \Big)_{\partial K}. 
  \end{equation}

Again, approximating $M$ with the lumped mass matrix  $M^{L}$, we obtain the following fully explicit scheme:
\begin{equation} \label{femod2n}
\begin{split}
 E^{k+1}_h =   & 
2  E^{k}_h  -  E^{k-1}_h  - \tau^2  (M^{L})^{-1}  G_1 E^{k}_h  -  \tau^2  (M^{L})^{-1}  G_2 E^{k}_h  \\
&+  
\tau^2 (M^{L})^{-1}  G_3 E^{k}_h + \tau^2  (M^{L})^{-1}  S_{\partial K} . 
\end{split}
\end{equation}
 
  \subsection{Finite difference formulation}
\label{sec:fdm}

 We recall now that from conditions (\ref{2.3}) it follows that in
 $\Omega_{FDM}$ the function $\varepsilon(x)= 1$.
This means that in
 $\Omega_{FDM}$ for the model problem \eqref{E_gauge1} the forward problem will be
\begin{align}\label{FDMmodel2}
\frac{\partial^2 E}{\partial t^2}  -  \Delta  E   = 0 &~~~ \mbox{in}~~
 \Omega_{FDM} \times (0,T), \\
  E(x,0) = f_0(x), ~~~E_t(x,0) = f_1(x) & ~~~\mbox{in}~~ \Omega_{FDM},  \\
   E = 0 &~~~\mbox{on}~ \partial \Omega \times (0,T),\\
\partial _{n} E = \partial _{n} E_{FEM} & ~~~\mbox{on}~ \partial \Omega_{FDM}  \times (0,T).
\end{align}

Using standard finite difference discretization of the
equation~(\ref{FDMmodel2}) in $\Omega_{FDM}$ we obtain the following explicit
scheme for the solution of the forward problem: 
\begin{equation}  \label{fdmschemeforward}  
  E_{l,j,m}^{k+1} =  \tau^2 \Delta E_{l,j,m}^k + 2 E_{l,j,m}^k - E_{l,j,m}^{k-1}.
\end{equation}
In the, system of, equations above,
$E_{l,j,m}^{k}$ is the solution at the time iteration $k$ at the
discrete point $(l,j,m)$,
 $\tau$ is the time step, and $ \Delta
E_{l,j,m}^k$ is the discrete Laplacian.

  Note that, in \eqref{fdmschemeforward}, the Dirichlet boundary consitions $E= E_{FEM}$  can be
  considered as well.

\section{Stability}
\label{sec:stability}
In this section we derive stability estimates for the semi-discrete approximations. For stability in 
$\Omega$ these estimates  are extensions of 
the stability approach derived for the continuous problem 
in \cite{BMaxwell}. As for the stability in $\Omega_{FEM}$ we get slightly different norms involving 
contributions corresponding to the 
reflexive  boundary: $\partial\Omega_{FEM}$.  
We use discrete version of a triple norm 
 induced by the weak variational formulation of 
 \eqref{eq1}, where we use the relation 
\eqref{eq6A1} (which is not necessary in the continuous case where 
$\nabla\cdot(\varepsilon E)=0$, however, in general $\nabla\cdot(\varepsilon E_h)\neq 0$): 

\emph{Find }$E\in {\textbf  W^E(\Omega)}$ \emph{ such  that}
\begin{equation}\label{eq1A1A}
  \begin{array}{l}
    \left ( \varepsilon \partial_{tt} E, {\textbf v} \right ) + (\nabla E,\nabla {\textbf v})+ 
    ((\nabla \varepsilon)\cdot  E, \nabla \cdot {\textbf v} ) 
    +((\varepsilon -1)\nabla \cdot E, \nabla \cdot {\textbf v})=0,
    \qquad \forall {\textbf v} \in {\textbf  W^E_0(\Omega)} \\
\\
\qquad   E(\cdot,0) = {f_0}(\cdot)\quad \mbox{ and } \,\,\, \partial_t E(\cdot,0) = {f_1}(\cdot) \quad 
   \mbox{ in } \Omega. 
  \end{array}
\end{equation}
\textbf{Remark}
  
In general, in non-divergent free case, the bilinear form induced by \eqref{eq1A1A} is not coercive. Further 
$H^1(\Omega)$-conforming finite element may result in spurious solutions. A remedy is through 
modifying the equation by adding a gauge constrain of Coulomb-type, see, e.g. \cite{Monk0} and \cite{div_cor}. 
This is supplied by the "zero"-term:  $\nabla\cdot(\varepsilon E)=0$, in \eqref{eq1}, which we 
add in the continuous variational formulation in \eqref{eq6}. This, however, is not necessarily true 
in the discrete forms, e.g. in \eqref{eq6B},  where most likely 
 $\nabla\cdot(\varepsilon E_h)\neq 0.$ 
Taking $ {\textbf v}=\partial_t E$  in \eqref{eq1A1A},  (we used the boundary condition $E=0$ on 
$\partial\Omega$), yields 
\begin{equation}\label{eq1A1B}
  \begin{split}
    \left ( \varepsilon \partial_{tt} E, \partial_t E\right ) & + (\nabla E,\nabla\partial_t E)+ 
    ((\nabla \varepsilon)\cdot  E, \nabla \cdot \partial_t E) \\
   & +((\varepsilon -1)\nabla \cdot E, \nabla \cdot \partial_t E) \equiv 0, 
  \end{split}
\end{equation}
which, due to the fact that  $\varepsilon$ is independent of $t$, 
can be rewritten as 
\begin{equation}\label{eq1A1C}
  \begin{split}
\frac 12 \frac d{dt}\Big(\varepsilon \partial_t E, \partial_t E\Big) &+ 
\frac 12 \frac d{dt}\Big( \nabla E, \nabla E\Big) 
+\Big((\nabla \varepsilon)\cdot  E, \nabla \cdot \partial_t E\Big)\\
&+\frac 12 \frac d{dt}\Big((\varepsilon -1)\nabla \cdot E, \nabla \cdot E\Big)
=0. 
 \end{split}
\end{equation}

\textbf{Proposition} \cite{BMaxwell}
  
Let $\Omega\subset{\mathbb R}^n,\, n=2,3$ be a bounded domain with piecewise linear boundary 
$\partial\Omega$. Then, the equation \eqref{eq1} has a unique solution $E\in H^2(\Omega_T)$. Further 
Let $f_1\in L_\varepsilon^2(\Omega)$ and 
$f_0\in H^1(\Omega)\cap H_{\varepsilon-1}^1(\Omega)$, then there is a constant  $C_\varepsilon^t=C(\vert\vert \varepsilon \vert\vert, t)$ such that, $\forall t\in (0, T]$,  the following 
stability estimate holds true 

\begin{equation}\label{eq1A1E}
  \begin{split}
\vert\vert\vert {E}\vert\vert\vert_{\varepsilon}^2& := 
 \vert\vert \partial_t E \vert\vert_\varepsilon^2(t)+
 \vert\vert \nabla E \vert\vert^2(t)
+ \vert\vert \nabla \cdot E  \vert\vert^2_{\varepsilon-1}(t) \\
&\le 
C_\varepsilon^t\Big(\vert\vert f_1 \vert\vert_\varepsilon^2  +
\vert\vert \nabla f_0 \vert\vert^2
 + \vert\vert f_0  \vert\vert^2+
 \vert\vert \nabla\cdot  f_0 \vert\vert_{\varepsilon-1}^2\Big).
 \end{split}
\end{equation}

\begin{proof} The estimate \eqref{eq1A1E} is proved in \cite{BMaxwell}, Theorem 4.1 by 
setting $s=1$ and $j\equiv 0$. 
 Integrating \eqref{eq1A1E} over the time interval $(0,t]$ we get the desired result. 
We omit the details. 
\end{proof}

Below we translate this stability to the semi-discrete problem.

\subsection{Stability estimate for the semi-discrete problem in $\Omega$}

\label{sec:stabOmega}

The stability for the semi-discrete problem in $\Omega$ is basically as in the continuous case above 
where all $E$:s are 
replaced by $E_h$ with some relevant assumptions in the discrete  data, viz. 

\textbf{Lemma}\label{StabOmega}

{\sl Assume that the interpolants of the data $f_0$ and $f_1$: $f_{0,h}$ and 
$f_{1,h}$ satisfy the regularity conditions $f_{1,h}\in L_\varepsilon^2(\Omega)$ and 
$f_{0,h} \in H^1(\Omega)\cap H_{\varepsilon-1}^1(\Omega)$, 
then for each $t\in (0,T]$}
\begin{equation}\label{StabAA}
\vert\vert\vert {E_h}\vert\vert\vert_{\varepsilon}^2(t)\le 
C_\varepsilon^t \Big( \vert\vert f_{1,h}  \vert\vert_\varepsilon^2
 +  \vert\vert \nabla f_{0,h} \vert\vert^2
 + \vert\vert f_{0,h} \vert\vert^2  +
 \vert\vert \nabla\cdot  f_{0,h} \vert\vert_{\varepsilon-1}^2\Big).
\end{equation}
where 
\begin{equation}\label{StabAB}
\vert\vert\vert {E_h}\vert\vert\vert_{ \varepsilon}^2(t):= 
 \vert\vert \partial_t E_h \vert\vert_\varepsilon^2(t)+
 \vert\vert \nabla E_h \vert\vert^2(t)
 + \vert\vert \nabla \cdot E_h \vert\vert^2_{\varepsilon -1}(t). 
\end{equation}

\subsection{Stability of the semi-discrete problem in $\Omega_{FEM}$ }

The stability of the semi-discrete problem in $\Omega_{FEM}$, 
relying on the variational formulation 
\eqref{eq6B}, and due to the appearance of the data function $g$,  is slightly 
different from \eqref{StabAA}.  We rewrite \eqref{eq6B}, in view of 
\eqref{eq6A1}, and with ${\textbf v}={\textbf v}_h=\partial_t E_h$ as: 
 {\sl Given $ E_h(\cdot,0)={f_0}_h$, $ \partial_t E_h(\cdot,0) ={f_1}_h$, and $g_h$, find 
$E_h\in  {\textbf W_h^E(\Omega_{FEM})}$ such that}
\begin{equation}\label{StabAC}
  \begin{split}
    \left ( \varepsilon \partial_{tt} E_{h}, \partial_t E_{h} \right ) &+
     (\nabla E_h,\nabla \partial_t E_{h})+ 
    ((\nabla \varepsilon) \cdot E_h, \nabla \cdot \partial_t E_{h})\\
    &+((\varepsilon -1)\nabla \cdot E_h, \nabla \cdot  \partial_t E_{h})  
= 
\langle g_{h}, \partial_t E_{h}\rangle_{\partial \Omega_{FEM}} . 
  \end{split}
\end{equation}
To deal with the $(\nabla\varepsilon)\cdot E_h$-term we rewrite 
\eqref{StabAC} in its original form as \eqref{eq6} for $\Omega_{FEM}$ and with ${\mathbf v}=\partial_t E_h$: 
\begin{equation}\label{StabAC1}
  \begin{split}
    \left ( \varepsilon \partial_{tt} E_{h}, \partial_t E_{h} \right ) &+
     (\nabla E_h,\nabla \partial_t E_{h})+ 
    (\nabla \cdot(\varepsilon  E_h), \nabla \cdot \partial_t E_{h})\\
    &-(\nabla \cdot E_h, \nabla \cdot  \partial_t E_{h})  
= 
\langle g_{h}, \partial_t E_{h}\rangle_{\partial \Omega_{FEM}}. 
 \end{split}
\end{equation}
Once again, in view of Theorem 4.1 in \cite{BMaxwell}, as in the case of stability in $\Omega$, we end up with 
the following time derivative form in $L_2(\Omega_{FEM})$ norms: 

\begin{equation}\label{StabAE}
  \begin{split}
\frac 12 \frac d{dt} & \Big(\vert\vert \partial_t E_h \vert\vert_\varepsilon^2+
\vert\vert \nabla E_h \vert\vert^2  +
 \vert\vert E_h \vert\vert^2+
\vert\vert \nabla \cdot E_h \vert\vert^2_{\varepsilon-1}\Big)+ \\
& \le \frac 12  \vert\vert g_h \vert\vert^2_{\partial\Omega_{FEM}}+
\frac 12 \vert\vert \partial_t E_h \vert\vert^2_{\partial\Omega_{FEM}}. 
 \end{split}
\end{equation}
Hence, integrating \eqref{StabAE} over the time interval $(0,t)$, we get $\forall t\in [0, T]$, 

\begin{equation}\label{StabAF}
  \begin{split}
\vert\vert\vert {E_h}\vert\vert\vert_{\varepsilon}^{2, FEM} (t):= & 
\Big( \vert\vert \partial_t E_h \vert\vert_\varepsilon^2+
\vert\vert \nabla E_h  \vert\vert^2
+ \vert\vert E_h \vert\vert^2 +
\vert\vert \nabla \cdot E_h  \vert\vert^2_{\varepsilon-1}\Big)(t)\\
\le &  \Big( \vert\vert \partial_t E_h \vert\vert_\varepsilon^2+
\vert\vert \nabla E_h \vert\vert^2 +
\vert\vert E_h \vert\vert^2+
\vert\vert \nabla \cdot E_h \vert\vert^2_{\varepsilon-1}\Big)(0)\\
&+ \int_0^t \vert\vert g_h \vert\vert^2_{\partial\Omega_{FEM}}+
 \int_0^t \vert\vert \partial_t E_h \vert\vert^2_{\partial\Omega_{FEM}}.
 \end{split}
\end{equation}

\textbf{Remark}

    We don't have electric conductivity: $\sigma$-term on the right hand side here. 
With the presence of $\sigma$ as in \eqref{eq1}, the associated assumptions are 
  $\sigma(x)\ge 1$ in $\Omega_{IN}$ and $\sigma(x)=0$ for 
$x\in \Omega_2 \cup \Omega_{OUT}$. As fot the boundary terms,  we may either use trace theorem and 
hide the $\partial\Omega_{FEM}$-terms in $ \vert\vert\vert {E_h}\vert\vert\vert_{\varepsilon}^{FEM} $ 
in \eqref{StabAF}, or redefine a modified version of 
$\vert\vert\vert {E_h}\vert\vert\vert_{\varepsilon}^{FEM} $ adding  terms corresponding to contributions 
from the boundary boundary. 
For the sake of generality we keep the two integrals as is and assume, 
for the boundary data, 
$\partial_t E_h\in L_2( \partial\Omega_{FEM})$.
Summing up we have the following stability estimate for the semi-discrete 
problem in $\Omega_{FEM}$: 

\textbf{Lemma}\label{SemiStab}

Under the following regularity assumptions on the interpolants for initial  conditions: 
$ f_{1,h}   \in L_\varepsilon^2(\Omega_{FEM})$, 
$f_{0,h} \in H^1(\Omega_{FEM})$, 
$\nabla\cdot   f_{0,h} \in L_{\varepsilon-1}^2(\Omega_{FEM})$, and with both boundary data: 
$g_h$, and $\partial_t E_h \in L_1\Big((0,T); L_2(\partial\Omega_{FEM})\Big)$, we have, 
for all $t\in (0, T]$, the following 
stability estimate for the semi-discrete $\Omega_{FEM}$ problem 
 \begin{equation}\label{StabAG}
  \begin{split}
\vert\vert\vert {E_h}\vert\vert\vert_{\varepsilon}^{2, FEM} (t)\le & 
\Big( \vert\vert f_{1,h} \vert\vert_\varepsilon^2
+ \vert\vert \nabla  f_{0,h} \vert\vert^2
 + \vert\vert f_{0,h} \vert\vert^2+
 \vert\vert \nabla\cdot f_{0,h} \vert\vert_{\varepsilon-1}^2\Big)(0)\\
 &+
  \int_0^t \vert\vert g_h \vert\vert^2_{\partial\Omega_{FEM}}+
  \int_0^t \vert\vert \partial_t E_h \vert\vert^2_{\partial\Omega_{FEM}}. 
 \end{split}
\end{equation}

 \textbf{Corollary}
 
We could write the right hand side in \eqref{StabAC} as 
$(g_{h}, \sqrt\varepsilon\partial_t E_{h}/\sqrt\varepsilon)_{\partial \Omega_{FEM}}$ . 
Then letting $C_{f_0, f_1}^2:=\Big( \vert\vert f_{1,h} \vert\vert_\varepsilon^2
+ \vert\vert \nabla f_{0,h}  \vert\vert^2
 + \vert\vert f_{0,h}   \vert\vert^2 +
 \vert\vert \nabla\cdot  f_{0,h}  \vert\vert_{\varepsilon-1}^2\Big)(0)$, the inequality 
 \eqref{StabAG} can be rewritten as 
\begin{equation}\label{StabAH}
\vert\vert\vert {E_h}\vert\vert\vert_{\varepsilon}^{2, FEM}(t) \le  
 C_{f_0, f_1}^2
 + \int_0^t \vert\vert g_h \vert\vert^2_{ \partial\Omega_{FEM}}+
  \int_0^t \vert\vert \frac{1}{\sqrt\varepsilon} \partial_t E_h \vert\vert^2_{\varepsilon, \partial\Omega_{FEM}}. 
\end{equation}
Thus by the definition of the triple norm and using Cauchy-Schwarz, Poincare and Gr\"{o}nwall's inequalities  
\begin{equation}\label{StabAHA}
\vert\vert E_h \vert\vert_{\varepsilon, \partial\Omega_{FEM}} \le \vert\vert \partial_t E_h \vert\vert_{\varepsilon, \partial\Omega_{FEM}}
\le C\Big(C_{f_0,f_1}+\int_0^T \vert\vert g_h \vert\vert_{ \partial\Omega_{FEM}}\Big) e^{T/\vert\vert \sqrt\varepsilon \vert\vert}. 
\end{equation} 
In a simlar way one may derive estimates of the gradient ($\nabla E_h $)-terms in the triple norm using 
the trace theorem, viz. 
\begin{equation}\label{StabAHB}
\begin{split}
\vert\vert \frac{1}{\sqrt\varepsilon} \partial_t E_h \vert\vert^2_{\varepsilon, \partial\Omega_{FEM}}&\le 
\sqrt[4]{8} \vert\vert \frac{1}{\sqrt\varepsilon} \partial_t E_h \vert\vert_{\varepsilon, \Omega_{FEM}}
\vert\vert \frac{1}{\sqrt\varepsilon} \partial_t E_h \vert\vert_{\varepsilon, W_2^1(\Omega_{FEM})} \\
&\le \frac{ \sqrt[4]{8}}2 \vert\vert \frac{1}{\sqrt\varepsilon} \partial_t E_h \vert\vert_{\varepsilon, \Omega_{FEM}}^2+
\frac {\sqrt[4]{8}}2 \vert\vert \frac{1}{\sqrt\varepsilon} \partial_t E_h \vert\vert_{\varepsilon, W_2^1(\Omega_{FEM})} ^2.
\end{split} 
\end{equation}
Now since both $\frac 12\sqrt[4]{8}<1$ likewise  $1/\sqrt \varepsilon <1 (\varepsilon >1)$ 
contributions from the right hand side  
terms can be hidden in 
corresponding terms of the triple norm, thus ending up with function and gradient terms estimates 
with  bounds depending on given parameters and functions 
$\sim {\mathcal M} (C_{f_0, f_1}, \, g_h,\, T,\, \varepsilon)$. 
We omit the  details. 
\section{Error estimates: Semi-Discrete (SD) problems}
\label{sec:errorest}
In what follows, and for future use in our model problems, 
we shall assume that $\partial_n E=-\partial_t E$ on $\partial \Omega$, 
which has the common value $g$ on $\partial\Omega_{FEM}$. 
Let now $\tilde{E}_h\in{\mathbf W}_h^{E}(\tilde\Omega)$, with 
$ \tilde\Omega =\Omega$ or  $ \tilde\Omega =\Omega_{FEM}$, be an spatial interpolant of the 
exact electric field $E$ and set 
\begin{equation}\label{SplitInterpol1}
e:=E-E_h=(E-\tilde{E}_h)+(\tilde{E}_h-E_h):=\eta+\xi.
\end{equation}
Then, assuming certain regularity of the data set, and with  
$\tilde{E}\in{\mathbf W}^{E}(\tilde\Omega)\cap H^s(\tilde\Omega)$, 
and with the spectral order $p\sim s \ge 1$, we can prove error estimates of the form 
\begin{equation}\label{ErrorEst1}
[\vert\vert e\vert\vert]_{\tilde\Omega}\le Ch^{p} \sim  Ch^s,\qquad \text{for}\quad 
 \tilde\Omega =\Omega\, \quad\mbox{or} \,\,\,  \tilde\Omega =\Omega_{FEM}. 
\end{equation}
  In this section, and to make a direct error estimate approach, without relying on the stability norm 
  defined in \cite{BMaxwell}, we use the equivialent norm $[\vert\vert \cdot \vert\vert]_{\tilde\Omega}$, slightly different 
  form the norm in \eqref{eq1A1E},  (see the term $ \vert\vert u \vert\vert^2_{  \vert \nabla\varepsilon \vert}$), and 
  directly obtained from the equation \eqref{eq6SEM}: 
  \begin{equation}\label{ErrorEst1Norm}
[\vert\vert u \vert\vert]_{\tilde\Omega}^2:= \vert\vert \partial_t u \vert\vert_\varepsilon^2+
\vert\vert \nabla u  \vert\vert^2 + \vert\vert u \vert\vert^2_{\vert \nabla\varepsilon \vert}
+ \vert\vert \nabla \cdot u \vert\vert^2_{\vert \nabla\varepsilon \vert +\varepsilon-1},
\quad
 \tilde\Omega =\Omega\, \,\,\,  \text{or}\quad  \tilde\Omega =\Omega_{FEM}. 
\end{equation}
Further, by the coercivity modification, see e.g.  \cite{Monk0}, there is a constant $C_\varepsilon$   such that 
 \begin{equation}\label{EquivCoerciv1}
[\vert\vert u \vert\vert]_{\tilde\Omega}^2\sim \vert\vert\vert u \vert\vert\vert_{\tilde\Omega}^2
\le  C_\varepsilon B_{\tilde\Omega}(u,u), 
\quad
 \tilde\Omega =\Omega\, \,\,\,  \text{or}\quad  \tilde\Omega =\Omega_{FEM}.
\end{equation}
Finally 
\begin{equation}\label{ErrorEst2Norm}
[\vert\vert u \vert\vert]_{\Omega_T}^2:= 
\int_0^T [\vert\vert u \vert\vert]_{\tilde\Omega}^2\, ds, 
\quad
 \Omega_T =\Omega\times [0,T]\, \,\,\,  \text{or}\quad  \Omega_T =\Omega_{FEM}\times[0,T], 
\end{equation}
likewise 
\begin{equation}\label{ErrorEst3Norm}
[\vert\vert\vert u \vert\vert]_{\Omega_T}^2:=\int_0^T[\vert\vert u \vert\vert]_{\tilde\Omega}^2,
\quad
 \Omega_T =\Omega\times [0,T]\, ds \,\,\,  \text{or}\quad  \Omega_T =\Omega_{FEM}\times[0,T].
\end{equation}

  
\textbf{Remark}
 
The original problem, with the presence of the electric conductivity terrm : $\sigma\partial_t E$\, 
($\sigma\neq 0$) on the right hand side, would behave as of parabolic type, 
(actually, quasi-parabolic, due to the presence of $\partial_{tt}E$-term).
Then in 
\eqref{ErrorEst1}, and for $E\in H^s(\tilde \Omega)$, $p\sim s$. But in our current consideration 
$\sigma\equiv 0$, and the problem is viewed as a system of wave equations 
(componenmtwise for E:s) and hence hyperbolic. On the other hand finite elements for 
 the scalar (non-system) hyperbolic problems has been considered in various studies by several 
 authors, showing that 
the best convergence one can hope is obtained  using, e.g. discontinuous 
 Galerkin (see \cite{Johnson_Pitkaranta}), which yields 
 \begin{equation}\label{ErrorEst1A}
[\vert\vert e\vert\vert]_{\tilde\Omega}\le Ch^{s-\theta},\qquad \text{for}\quad 
 \tilde\Omega =\Omega\, \quad\mbox{or} \,\,\, \Omega_{FEM}, 
\end{equation}
instead of 
   \eqref{ErrorEst1} and with $\theta=1/2$, whereas the finite difference approach for the same, hyperbolic type, problem 
   is more accurate and satisfies  \eqref{ErrorEst1}.  

Below we use the very similar argument to 
derive \eqref{ErrorEst1} for the spatial domains $\Omega$ and $\Omega_{FEM}$. 
   


\subsection{Error estimates: SD problem in $\Omega$}
\label{errerOmega} 

\textbf{Theorem}

\label{ErrorEstThm1} For $E\in {\mathbf W}^2(\Omega)$ and 
continuous  piecewise polynomial approximation, 
assuming  
$$
\vert\vert f_1 \vert\vert_{L^2_\varepsilon  (\Omega)}^2+ \vert\vert f_0 \vert\vert_{H^1(\Omega)}^2+
\vert\vert f_0 \vert\vert_{\vert \nabla\varepsilon   \vert +\varepsilon-1}^2(\Omega)+
\vert\vert f_0 \vert\vert_{H^1_{\varepsilon-1}(\Omega)}^2\le C,
$$
then there is a constant $C$ such that 
\begin{equation}\label{ErrorEs2}
[\vert\vert e\vert\vert]\le Ch^{}, 
\end{equation}

\begin{proof}
We start with the straightforward estimate for $\eta$ in \eqref{SplitInterpol1},  using  
interpolation error:
\begin{equation}
\left\| u-u_h^\mathrm{I}\right\| _{L_{2}\left( \Omega \right) }\leq C_%
\mathrm{I} \left\|h~ \nabla u\right\| _{L_{2}\left( \Omega \right) }.
\label{2.6A}
\end{equation}
Note that if $u\in C^2(\Omega)$, then  the continuous interpolation, \eqref{2.6A} 
 is improved, and  
 
 \begin{equation*}
\left\| u-u_h^\mathrm{I}\right\| _{L_{2}\left( \Omega \right) }
\leq C_\mathrm{I}  h^2  \vert\vert D_x^2 u \vert\vert. 
\end{equation*}
However, such improvement can not survive, e.g.  in  approximating with discontinuous interpolation 
where jump terms   $\left [\frac{\partial u_h}{\partial n} \right ]$ are introduced in the outward normal directions to elements, 
and 
$$
  D_x^2 u \leq \frac{\left [\frac{\partial u_h}{\partial n} \right ] }{h}.
$$
Consequently,  returning to   \eqref{2.6A}  we get:
\begin{equation*}
  \left\| u-u_h^\mathrm{I}\right\| _{L_{2}\left( \Omega \right) }\leq C_\mathrm{I}
  h  \left |  \left [\frac{\partial u_h}{\partial n} \right ]  \right|.
\end{equation*}
Hence, we have the following estimate for the interpolation error: 
\begin{equation}\label{Interpol1}
\begin{split}
[\vert\vert \eta\vert\vert]_{\varepsilon}^2 (t)=& 
[\vert\vert E  - \tilde{E}_h\vert\vert]_{\varepsilon}^2= 
\vert\vert \partial_t(E-\tilde{E}_h) \vert\vert_\varepsilon^2+  \vert\vert \nabla(E-\tilde{E}_h)\vert\vert^2
\\
&+
\vert\vert E-\tilde{E}_h \vert\vert^2_{\vert \nabla\varepsilon \vert}
+ \vert\vert \nabla\cdot(E-\tilde{E}_h) \vert\vert^2_{(\vert \nabla\varepsilon \vert + \varepsilon-1)}\\
\le &
 h^2\Big ( \vert\vert \partial_t E  \vert\vert^2_{\varepsilon}+
\vert\vert \nabla E  \vert\vert^2 + \vert\vert E \vert\vert^2_{\vert \nabla\varepsilon \vert}+ \vert\vert \nabla \cdot E \vert\vert^2_{\vert \nabla\varepsilon \vert +\varepsilon-1}\Big)(t) 
\le  (C^t_\varepsilon)^2 h^{2},
\end{split}
\end{equation}
where the last two inequalities are just the consequences of the interpolation error and regularity of the exact solution, respectively. So that we can deduce that the interpolation error is:
\begin{equation}\label{Interpol2}
[\vert\vert \eta\vert\vert]_{\varepsilon}(c)\sim C^t _\varepsilon h.
\end{equation}
To proceed further we assume continuous variational formulation, i.e.   
the continuous version of \eqref{eq6SEM}:
$$
0 = B_\Omega(E, {\mathbf v})=0 \qquad \forall \,{\mathbf v}\in {\mathbf W}_h^{E}(\Omega). 
$$
Hence we can write  for $\xi = \tilde{E}_h-E_h$
\begin{equation}\label{ErrorEs3}
B_\Omega(\xi, \xi)= B_\Omega(\tilde E_h-E_h, \xi)=
B_\Omega(\tilde E_h, \xi)
=B_\Omega(E-\tilde E_h,  \xi)=B_\Omega(\eta, \xi)
\end{equation}
and hence, in a time interval $(0,t)\subset (0,T)$ we have 
\begin{equation}\label{ErrorEs4} 
[\vert\vert \xi\vert\vert]_{\varepsilon}^2 =
[\vert\vert \tilde E_h-E_h \vert\vert]_{\varepsilon}^2\le C_\varepsilon^t 
\int_0^TB_\Omega(\xi, \xi) \, dt =C_\varepsilon^t
\int_0^TB_\Omega(\eta, \xi) \, dt, 
\end{equation}
where 
\begin{equation}\label{ErrorEs5} 
\begin{split}
B_\Omega(\eta, \xi) =& (\varepsilon \eta_{tt}, \xi) + 
(\nabla \eta, \nabla\xi)+((\nabla\varepsilon)\cdot\eta,\nabla\cdot\xi) \\
&+\Big((\varepsilon -1)\nabla\cdot\eta,\nabla\cdot\xi\Big)
:=\sum_{k=1}^4J _k(t).
\end{split}
\end{equation}
We estimate each $\int_0^t J_k(s)\, ds, $ for $k=1,2, 3, 4$, separately. As for $J_1$, partial integration, with zero boundary condition, yields 
\begin{equation}\label{ErrorEs6} 
\begin{split}
\int_0^tJ_1(s)\, ds &=\int_0^t(\varepsilon \eta_{ss}, \xi) \, ds=
-\int_0^t(\sqrt\varepsilon \eta_{s}, \sqrt\varepsilon \xi_s)\, ds \\
&\le \int_0^t \vert\vert \eta_{s}  \vert\vert^2_\varepsilon \, ds+\frac 14 \int_0^t \vert\vert \xi_s \vert\vert^2_\varepsilon \, ds. 
\end{split}
\end{equation}
Direct estimates for $J_2$, $J_3$ (with some formal manipulations), and $J_4$-terms give 
\begin{equation}\label{ErrorEs7} 
\int_0^t J_2(s)\, ds =\int_0^t (\nabla \eta, \nabla\xi)\le 
\int_0^t \vert\vert \nabla \eta \vert\vert^2\, ds+\frac 14\int_0^t \vert\vert \nabla\xi \vert\vert^2\, ds,
\end{equation}

\begin{equation}\label{ErrorEs8} 
\int_0^t J_3(s)\, ds =\int_0^t ((\nabla \varepsilon)\cdot\eta, \nabla\cdot\xi)\le 
\int_0^t \vert\vert \eta \vert\vert^2_{\vert \nabla\varepsilon \vert}\, ds+\frac 14\int_0^t \vert\vert \nabla\cdot\xi  \vert\vert^2_{\vert \nabla\varepsilon \vert}\, ds,
\end{equation}

\begin{equation}\label{ErrorEs9} 
\begin{split}
\int_0^tJ_4(s)\, ds &=\int_0^t \Big((\varepsilon -1) \nabla\cdot\eta, \nabla\cdot\xi\Big)\, ds=
\int_0^t(\sqrt{\varepsilon -1}\nabla\cdot\eta, \sqrt{\varepsilon-1} \nabla\cdot\xi)\, ds \\
&\le \int_0^t \vert\vert  \nabla\cdot\eta \vert\vert^2_{\varepsilon -1}\, ds
+\frac 14 \int_0^t \vert\vert \nabla\cdot\xi \vert\vert^2_{\varepsilon -1}\, ds. 
\end{split}
\end{equation}
Thus by a kick back argument 
all $\xi$-norms on the right hand side, 
can be hidden in the corresponding  terms in $[\vert\vert{\xi}\vert\vert]_{\varepsilon}^2$ leading to a 
$\eta$ estimate for $\xi$: 
\begin{equation}\label{ErrorEs12t} 
[\vert\vert{\xi}\vert\vert]_{\varepsilon}(t)\le C_\varepsilon^t [\vert\vert{\eta}\vert\vert]_{\varepsilon}(t).
\end{equation} 
Now, recalling 
\eqref{Interpol2} we get the desired result. 
\end{proof} 

\textbf{Proposition}

With the same assumption as in the theorem \ref{ErrorEstThm1} above we have the convergence rate of the time derivative $e_t$ for the error:

\begin{equation}\label{ErrorEs12tB} 
[\vert\vert{e_t}\vert\vert]_{\varepsilon}(t)\le C_\varepsilon^t  h.
\end{equation} 

\begin{proof}
Evidently the interpolation estimates in the proof of theorem \ref{ErrorEstThm1} also yield for $\eta_t$, and 
\begin{equation}\label{Eta_t1}
[\vert\vert{\eta_t}\vert\vert]_{\varepsilon}(t)\le C_\varepsilon^t h.
\end{equation} 
Note that we have no time discretization yet. The remaining step is to show that 
\begin{equation}\label{ErrorEs12tBA} 
[\vert\vert{\xi_t}\vert\vert]_{\varepsilon}(t)\le C_\varepsilon^t [\vert\vert{\eta_t}\vert\vert]_{\varepsilon}(t).
\end{equation} 
The same proceedure with $\eta$ and $\xi$ replaced by $\eta_t$ and $\xi_t$, respectively, yields
\begin{equation}\label{ErrorEs3t}
\begin{split}
B_\Omega(\xi_t, \partial_t\xi)&= B_\Omega((\tilde E_h-E_h)_t, \partial_t\xi)=
B_\Omega(\tilde E_{h,t}, \partial_t\xi)\\
&=B_\Omega((E-\tilde E_h)_t,  \partial_t\xi)=B_\Omega(\eta_t, \partial_t\xi)
\end{split}
\end{equation}
and hence 
\begin{equation}\label{ErrorEs4t} 
[\vert\vert \xi_t\vert\vert]_{\varepsilon}^2 =
[\vert\vert\partial_t (\tilde E_{h}- E_h) \vert\vert]_{\varepsilon}^2\le 
C_\varepsilon^t \int_0^TB_\Omega(\xi_t, \xi_t) \, dt =
C_\varepsilon^t \int_0^TB_\Omega(\eta_t, \xi_t) \, dt, 
\end{equation} 
where 
\begin{equation}\label{ErrorEs5t} 
\begin{split}
B_\Omega(\eta_t, \partial_t\xi) =& (\varepsilon \eta_{ttt}, \xi_t) + 
(\nabla \eta_t, \nabla\xi_t)+((\nabla\varepsilon)\cdot\eta_t,\nabla\cdot\xi_t) \\
&+\Big((\varepsilon -1)\nabla\cdot\eta_t,\nabla\cdot\xi_t\Big)
:=\sum_{k=1}^4 I_k(t).
\end{split}
\end{equation}
Mimiking the above procedure we estimate $\int_0^t I_k(s)\, ds$-terms  for $k=1,2, 3, 4$, viz. 
\begin{equation}\label{ErrorEs6t} 
\begin{split}
\int_0^t I_1(s)\, ds &=\int_0^t(\varepsilon \eta_{sss}, \xi_s) \, ds=
-\int_0^t(\sqrt\varepsilon \eta_{ss}, \sqrt\varepsilon \xi_{ss})\, ds \\
&+\varepsilon\eta_{ss}(t)\xi_s(t)-\varepsilon\eta_{ss}(0)\xi_s(0)\le \int_0^t \vert\vert  \eta_{ss} \vert\vert ^2_\varepsilon \, ds+\frac 14 \int_0^t \vert\vert \xi_{ss} \vert\vert ^2_\varepsilon \, ds, 
\end{split}
\end{equation}
where, with the continuous in time, the spatial discrete errors for $\eta_{ss}$ and $\xi_s$ 
are assumed to be zero for all  $t\in [0,T]$. 

\begin{equation}\label{ErrorEs7t} 
\int_0^t I_2(s)\, ds =\int_0^t (\nabla \eta_s, \nabla\xi_s)\le 
\int_0^t \vert\vert \nabla \eta_s \vert\vert ^2\, ds+\frac 14\int_0^t \vert\vert \nabla\xi_s \vert\vert^2\, ds,
\end{equation}

\begin{equation}\label{ErrorEs8t} 
\int_0^t I_3(s)\, ds =\int_0^t ((\nabla\varepsilon) \cdot\eta_s, \nabla\cdot\xi_s)\le 
\int_0^t \vert\vert \eta_s \vert\vert ^2_{\vert \nabla\varepsilon \vert}\, ds
+ \frac 14\int_0^t \vert\vert \nabla\cdot\xi_s \vert\vert^2_{\vert \nabla\varepsilon \vert}\ \, ds,
\end{equation}

\begin{equation}\label{ErrorEs9t} 
\begin{split}
\int_0^tI_4(s)\, ds &=\int_0^t \Big((\varepsilon -1) \nabla\cdot\eta_s, \nabla\cdot\xi_s\Big)\, ds=
\int_0^t(\sqrt{\varepsilon -1}\nabla\cdot\eta_s, \sqrt{\varepsilon-1} \nabla\cdot\xi_s)\, ds \\
&\le \int_0^t \vert\vert \nabla\cdot\eta_s \vert\vert^2_{\varepsilon -1}\, ds
+\frac 14 \int_0^t \vert\vert  \nabla\cdot\xi_s \vert\vert^2_{\varepsilon -1}\, ds. 
\end{split}
\end{equation}
Note that due to the vanishing boundary condition, the contribution from the boundary is not present. 
This however can be inserted by considering a modified {\sl triple} norm including, e.g., reflecting 
boundaries as in the case of $\Omega_{FEM}$ below. 
Hence, once again, a kick-back argument, with all $\xi_s$ and $\xi_{ss}$ weighted-norms on the left hand side 
are hidden in the corresponding  terms in $[\vert\vert{\xi_t}\vert\vert]_{\varepsilon}^2$ giving the 
$\eta_t$ estimate \eqref {ErrorEs12tBA} for $\xi_t$,
which, combined with 
\eqref{Eta_t1}, gives the desired result. 
\end{proof}

\subsection{Error estimates: SD problem in $\Omega_{FEM}$} \label{errorOmegaFEM}

The estimates here are mostly the same as  those of ther previous subsection. However, 
here we have a reflexive boundary condition on $\partial\Omega_{FEM}$. Hence, the 
estimates contain an extra contribution from the boundary 
( in the previous subsection, we have only considered 
the zero boundary condition for the whole $\Omega$). Here, we include the procedure containg 
boudarry trerm estimates, which can be mimiked in the case of the refexive boundary condition in whole $\Omega$. 

\textbf{Theorem}\label{ErroreEstThm1}

Let $E\in {\mathbf W}^2(\Omega_{FEM})$ and consider 
the continuous piecewise polynomial approximation for the solution of the problem 
\eqref{eq1fem}. Furthermore, assume that $f_1\in L_{2, \varepsilon}(\Omega_{FEM})$, 
$f_0\in H^1(\Omega_{FEM})\cap H^1_{\nabla\varepsilon}(\Omega_{FEM})\cap H^1_{\varepsilon-1}(\Omega_{FEM})$, 
and 
$g\in L_2(\partial\Omega_{FEM})$.  
Then there is a constant $C$ such that 
\begin{equation}\label{ErrorEs2}
[\vert\vert e\vert\vert]_{\varepsilon, \Omega_{FEM}}\le C_\varepsilon^t h^{}. 
\end{equation}

\begin{proof}
Following the same procedure as the error estimates in $\Omega$, and letting now 
\begin{equation}\label{ErrEstFEM1}
(E-E_h)_{\Omega_{FEM}}:=
(E-\tilde E_h)_{\Omega_{FEM}}+(\tilde E_h-E_h)_{\Omega_{FEM}}
:=\rho+\theta, 
\end{equation} 
we need to estimate a triple norm of the form 
\begin{equation}\label{Interpol_FEM}
\begin{split}
&[\vert\vert E-E_h \vert\vert]_{\varepsilon, \Omega_{FEM}}^2:= 
\vert\vert \partial_t(E-{E}_h) \vert\vert_{\varepsilon, \Omega_{FEM}}^2+ 
\vert\vert \nabla(E-{E}_h) \vert\vert_{\Omega_{FEM}}^2\\
\qquad &+ 
\vert\vert E - \tilde{E}_h \vert\vert^2_{\vert \nabla\varepsilon \vert, \Omega_{FEM}} 
+ \vert\vert \nabla\cdot(E- {E}_h) \vert\vert^2_{\vert \nabla\varepsilon \vert + \varepsilon-1, \Omega_{FEM}}\\
&+
\int_0^t \vert \partial_t(E-{E}_h)\vert^2_{\partial\Omega_{FEM}}(s)\, ds 
+\int_0^t \vert g-g_h \vert^2_{\partial\Omega_{FEM}}(s)\, ds.
\end{split}
\end{equation}
Assuming that, at the boundary $\partial\Omega$,  $g$ is as regular as $E$,  
the linear interpolation error reads: 
\begin{equation}\label{Interpol_FEM2}
\vert\vert\vert \rho\vert\vert\vert_{\varepsilon, \Omega_{FEM}}(t)\sim C_\varepsilon^t h_{\Omega_{FEM}}, 
\end{equation}
Then following the same procedure as above  we have for both \eqref{eq6B}
and its continuous version, and  with $\zeta=E$ and $ \zeta=E_h$ corresponding to 
$G=g,$ and  $G=g_h$, respectively, we have 
\begin{equation}\label{eq6B_FEM}
  \begin{split}
    B_{\Omega_{FEM}}(\zeta, {\textbf v}):&=\left ( \varepsilon \partial_{tt} \zeta, {\textbf v} \right ) + 
    (\nabla \zeta,\nabla {\textbf v})+ 
    ((\nabla \varepsilon )\cdot\zeta), \nabla \cdot {\textbf v} ) 
    -((\varepsilon-1)\nabla \cdot \zeta, \nabla \cdot {\textbf v}) \\
    & =\langle G, {\textbf v}\rangle_{\partial \Omega_{FEM}} ,  \quad 
      \forall  {\textbf v} \in  {\textbf W_h^E(\Omega_{FEM})}.
    \end{split} 
 \end{equation}
with the associated data. 
Hence we can write 
\begin{equation}\label{ErrorEs3_FEM}
\begin{split}
B_{\Omega_{FEM}}(\theta, \theta)&= B_{\Omega_{FEM}}(\tilde E_h-E_h, \theta) \\
 &-B_{\Omega_{FEM}}(E-\tilde E_h,  \theta)+B_{\Omega_{FEM}}(E-\tilde E_h,  \theta)\\
&=-B_{\Omega_{FEM}}(e, \theta)+B_{\Omega_{FEM}}(\rho, \theta) \\
&=
\langle g_h-g,\theta\rangle_{\partial\Omega_{FEM}}+B_{\Omega_{FEM}}(\rho, \theta). 
\end{split}
\end{equation}  
In the triple norm form this yields 
\begin{equation}\label{ErrorEs4_FEM} 
\begin{split}
[\vert\vert \theta\vert\vert]_{\varepsilon, \Omega_{FEM}}^2(t) &=
[\vert\vert \tilde E_h-E_h \vert\vert]_{\varepsilon, \Omega_{FEM}}^2(t)\le C_\varepsilon^t 
\int_0^t B_{\Omega_{FEM}}(\theta, \theta) \, ds \\
&=
\int_0^t \langle g_h -g, \theta\rangle_{\partial\Omega_{FEM}} ds+\int_0^tB_{\Omega_{FEM}}(\rho, \theta) \, ds\\
&\le 
\int_0^t \vert g-g_h \vert^2_{\partial\Omega_{FEM}}\, ds
+\frac 14  \int_0^t  \vert \theta \vert^2_{\partial\Omega_{FEM}}\, ds \\
&+  [\vert\vert \rho\vert\vert]_{\varepsilon, \Omega_{FEM}}^2(t) 
+\frac 14  [\vert\vert \theta\vert\vert]_{\varepsilon, \Omega_{FEM}}^2(t). 
\end{split} 
 \end{equation}
Now using the same procedure as in the proof of the previous theorem 
(with $\xi$ and $\eta$ replaced by $\rho$ and $\theta$) to bound all the involved norms 
and hiding both $\theta$-terms on the right, inside 
$[\vert\vert \theta\vert\vert]_{\varepsilon, \Omega_{FEM}}^2(t) $,  on the left hand side, 
together with estimates \eqref{Interpol_FEM2}
for $\rho$, and further assuming  corresponding estimates for 
$(g_h -g, \theta)_{\partial\Omega_{FEM}} $, we obtain (omitting  some details) 
that 
\begin{equation}\label{ErrorEs12t_FEM} 
[\vert\vert{\theta}\vert\vert]_{\varepsilon}^2(t)\le (C^t_\varepsilon)^2 
\Big([\vert\vert{\rho}\vert\vert]_{\varepsilon}^2(t)+
\int_0^t \vert \langle g_h -g, \theta\rangle \vert_{\partial\Omega_{FEM}} ^2\Big)\, ds. 
\end{equation} 
Now 
\eqref{Interpol_FEM2} and \eqref{ErrorEs12t_FEM}  together with the corresponding 
estimates for $\langle g_h -g, \theta\rangle_{\partial\Omega_{FEM}} $, give the desired result and completes the proof.  
\end{proof}

{
\begin{algorithm}[hbt!]
  \centering
  \caption{Domain decomposition algorithm for solution of Maxwell's equations \eqref{E_gauge1}. At every time step $k$ are performed the following operations:} 
  \vskip 0.3cm 
  \begin{algorithmic}[1]

    \STATE Compute $E^{k+1}$ in $\Omega_{FDM}$ using the explicit
  finite difference  scheme \eqref{fdmschemeforward} with known $E^k$, and
    $E^{k-1}$-values.

 \STATE Compute $E^{k+1}$ in $\Omega_{\rm FEM}$ by using the finite
 element scheme \eqref{femod2} with known $E^k, E^{k-1}$.
    
   \STATE For the finite
  difference method in $\Omega_{\rm FDM}$, use the values of the function $E^{k+1}$ at nodes
  $\omega_{\diamond}$ (green boundary of Figure \ref{fig:F1}) , which are computed using the finite element
  scheme \eqref{femod2n}, as a boundary condition at the inner boundary of $\Omega_{\rm FDM}$. 

  \STATE Apply appropriate boundary condition at the outer boundary of $\Omega_{\rm FDM}$.
  
    \STATE  For the finite element
  method in $\Omega_{\rm FEM}$,  use the values of the functions $E^{k+1}$ at
  nodes $\omega_{\rm o}$  (blue boundary of the Figure \ref{fig:F1} ), 
  which are computed using the finite difference
  scheme \eqref{fdmschemeforward} as a boundary condition.

   \STATE   Apply swap of the solutions for the computed function $E^{k+1}$
   to be able to perform the algorithm on a new time level $k$.
   
  \end{algorithmic}
\end{algorithm}
}

\section{Numerical  examples}

In this section we present numerical examples justifying theoretical results of the previous two sections. 
For convergence tests the domain decomposition algorithm (see Algorithm 2), 
implemented in the software package WavES \cite{waves}, was used.
We note that because of
using explicit FE and FD schemes in $\Omega_{\rm FEM}$ and $\Omega_{\rm FDM}$, correspondingly,  we need to choose time step $\tau$   according to the CFL stability condition \eqref{CFL1} derived in \cite{BR1} so that
the whole hybrid scheme remains stable.

Numerical tests are performed in time interval $(0,T)=(0,0.25)$ and in the  spatial  dimensionless
computational domain
 \begin{equation}\label{ABC}
 \Omega = \left\{ (x,y): x \in [0, 1], y  \in [0,1] \right\},
 \end{equation} 
  which is split into the finite element domain
 \begin{equation} \label{ABCD}
\Omega_{\rm FEM} = \left\{ (x,y): x \in [0.25,0.75], y \in [0.25, 0.75]  \right\}
 \end{equation}
and the finite difference domain $\Omega_{\rm FDM}$, thus 
$\Omega = \Omega_{\rm FEM} \cup \Omega_{\rm FDM}$, see Figure \ref{fig:F1}.

The   model problem in all our tests that is stated for the electric field  $E = (E_1, E_2)$  is as follows:
\begin{equation}\label{model}
  \begin{array}{ll}
    \varepsilon \partial_{tt} E +  \nabla (\nabla \cdot E) - \triangle E
    - \nabla \nabla \cdot (\varepsilon  E)  = F & \mbox{ in } \Omega \times (0, T), \\
    E(\cdot,0) =  0 \mbox{ and } \partial_t E(\cdot,0)
    =  0& \mbox{ in } \Omega, \\
    E =  0 & \mbox{ on } \partial \Omega \times (0,T).
  \end{array}
\end{equation}

We have the functions 
\begin{equation}\label{eq1exact}
  \begin{split}
    E_1 &=    \frac{1}{\varepsilon} 2 \pi \sin^2 \pi x  \cos \pi y  \sin \pi y ~ \frac{t^2}{2}, \\
   E_2 &=  - \frac{1}{\varepsilon}   2 \pi   \sin^2 \pi y \cos \pi x   \sin \pi x   ~\frac{t^2}{2}, 
  \end{split}
\end{equation}
as the exact solution $E = (E_1, E_2)$ of the model problem
\eqref{model} with the source data $F=(F_1,F_2)$ which corresponds to
this exact solution.

\begin{figure}[h!]
\begin{center}
\begin{tabular}{cccc}
  {\includegraphics[scale=0.1, clip=]{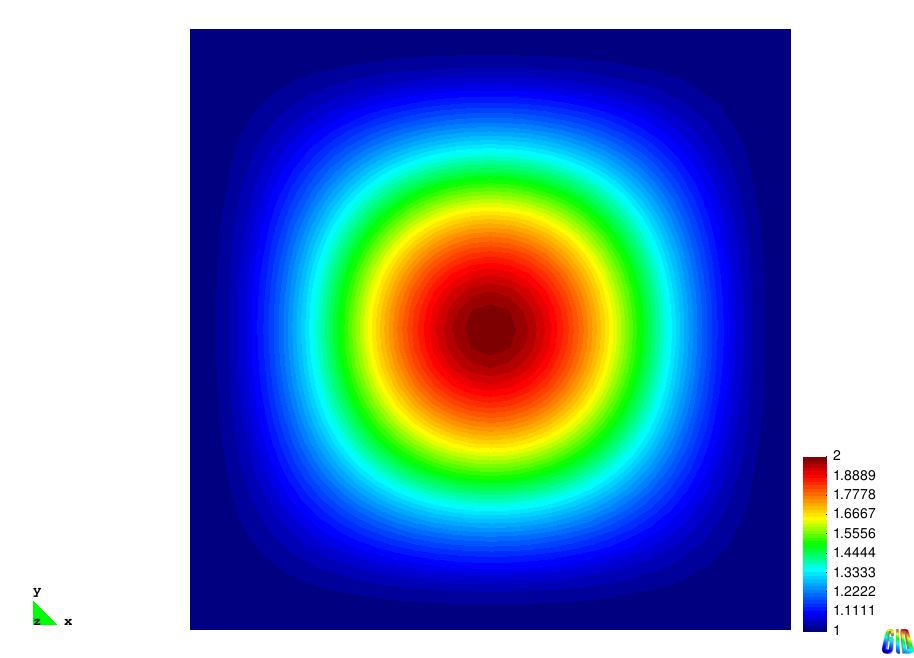}} &
  {\includegraphics[scale=0.1, clip=]{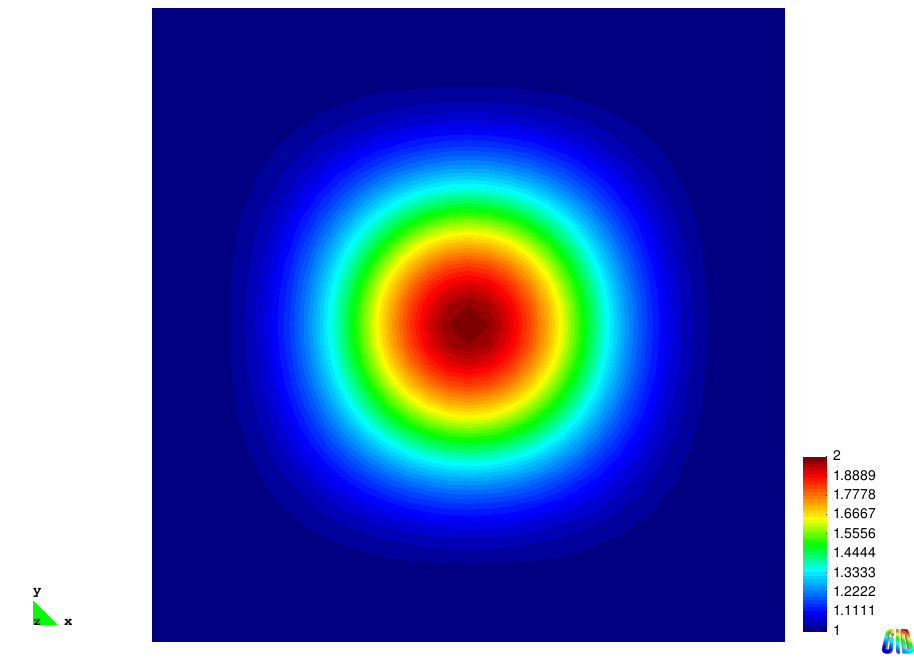}} &
  {\includegraphics[scale=0.1, clip=]{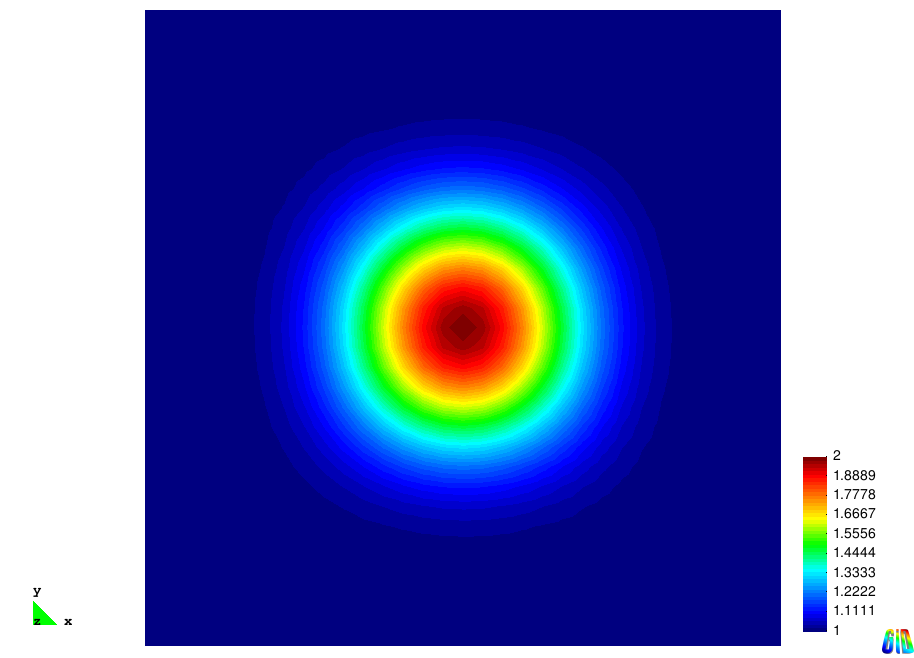}} &
   {\includegraphics[scale=0.1, clip=]{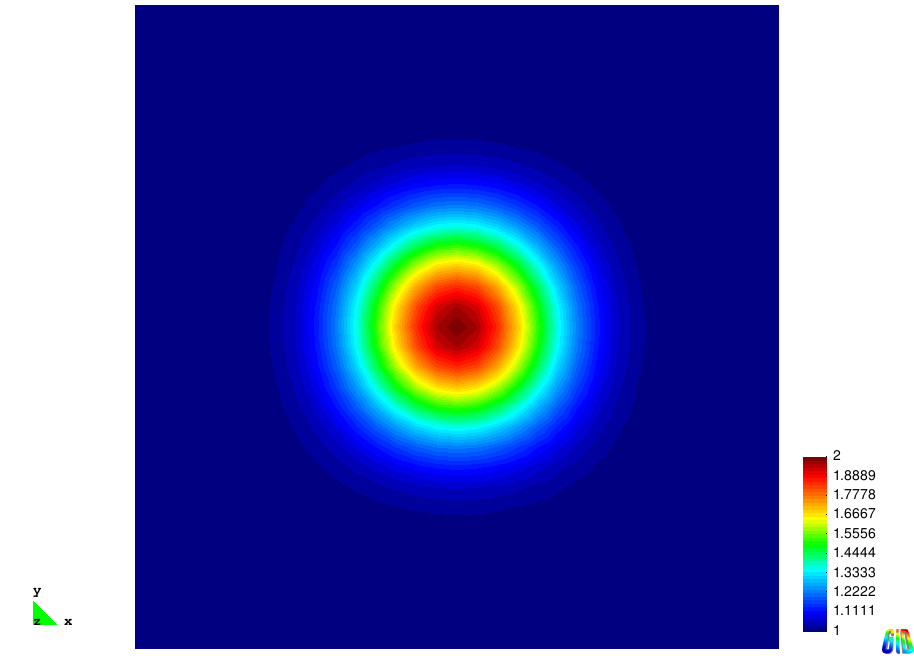}} \\
   $m=2$&  $m=4$ &  $m=6$ &   $m=8$
\end{tabular}
\end{center}
\caption{Function $\varepsilon(x,y)$ in the domain $\Omega_{\rm FEM}$ for different values of $m$ in \eqref{eps}.}
\label{fig:F5}
\end{figure}

 The function $\varepsilon$ in \eqref{model} is defined as
\begin{equation}\label{eps}
  \varepsilon(x,y)= \left \{
  \begin{array}{ll}
    1 + \sin^m \pi (2x-0.5) \cdot \sin^m \pi (2y-0.5) & \textrm{in $\Omega_{FEM}$}, \\
    1 &  \textrm{in $\Omega_{FDM}$}.
  \end{array}
  \right.
\end{equation}
We choose $m = 2, 4, 6, 8$ in our numerical examples, see  Figure \ref{fig:F5}, for these functions
in the domain $\Omega_{\rm FEM}$. We note that the exact solution \eqref{eq1exact}  satisfies the divergence free condition $\nabla \cdot (\varepsilon E) =0$ for
  $\varepsilon$  defined by \eqref{eps}, the homogeneous
initial conditions, as well as the homogeneous Dirichlet conditions  for all times.

\begin{figure}[h!]
\begin{center}
\begin{tabular}{cc}
  {\includegraphics[scale=0.4, clip=]{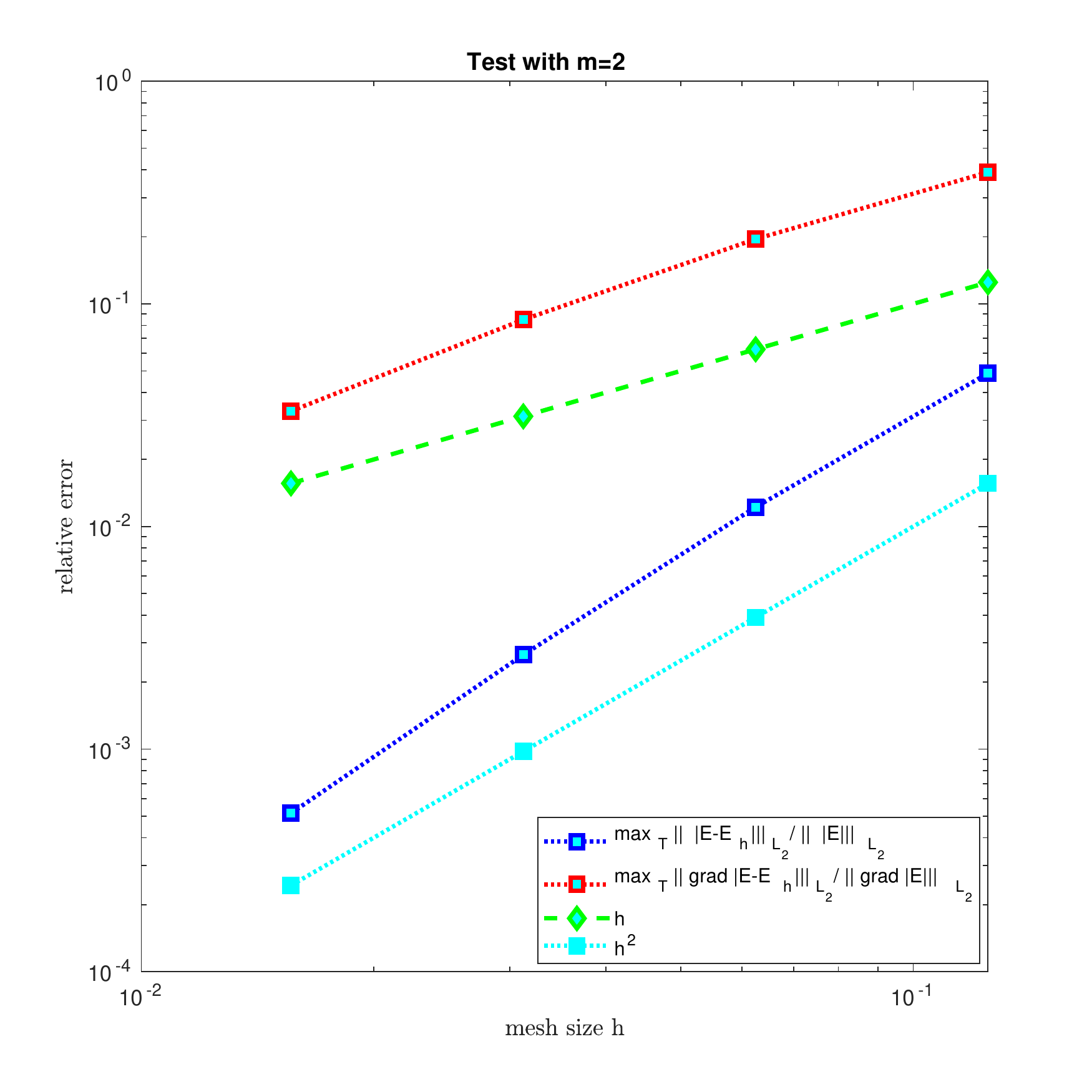}} &
  {\includegraphics[scale=0.4, clip=]{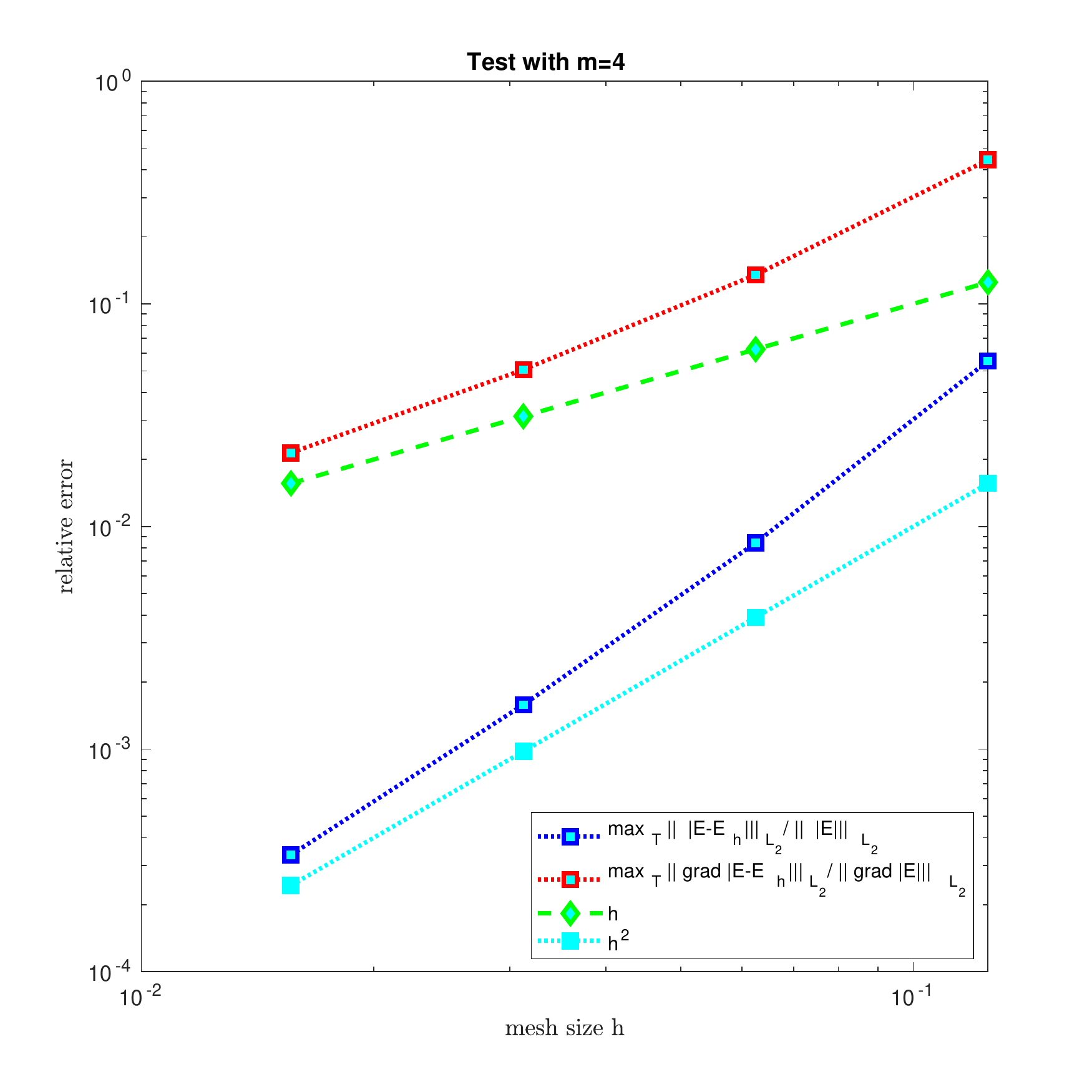}}   \\
  m=2 & m= 4 \\
   {\includegraphics[scale=0.4, clip=]{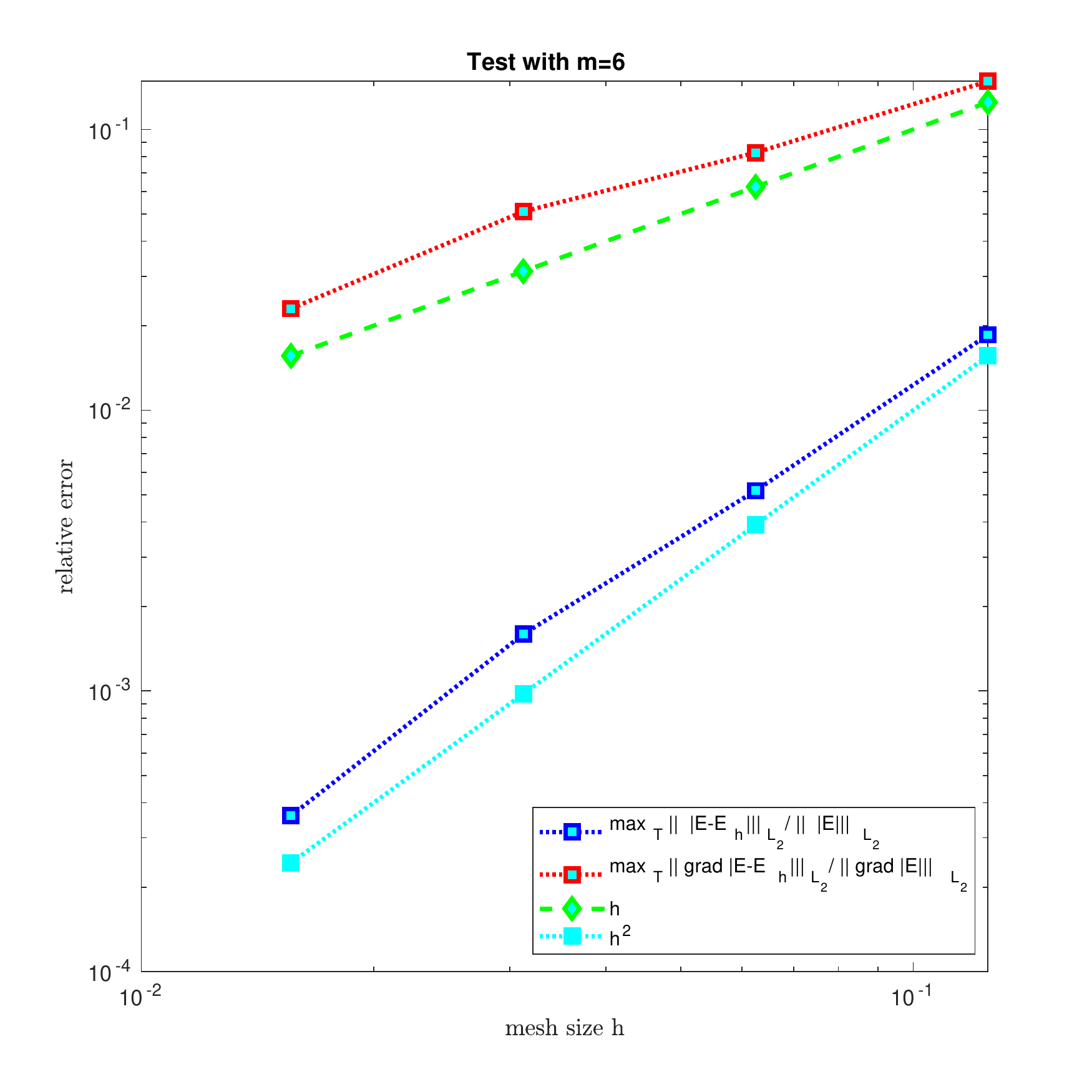}} &
   {\includegraphics[scale=0.4, clip=]{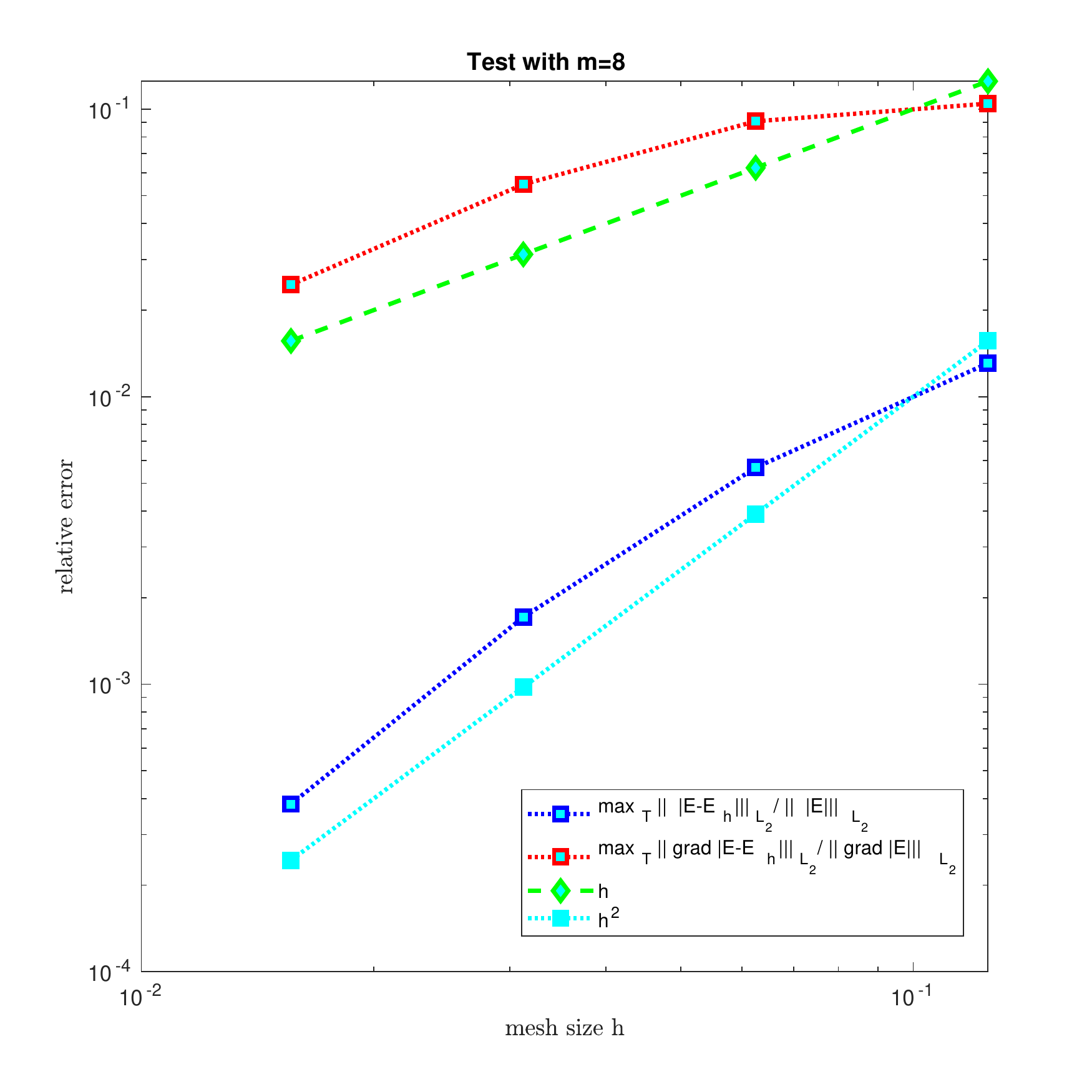}}   \\
   m=6 & m=8\\
\end{tabular}
\end{center}
\caption{Convergence  of the   relative $L_2$ and $H^1$ norms computed  via \eqref{norms}.}
\label{fig:F3}
\end{figure}

The computational domain $\Omega_{FEM} \times (0,T)$ was discretized
into triangular elements with mesh sizes $h_l= 2^{-l}, l=3,4,5,6$, and the mesh in
$\Omega_{FDM} \times (0,T)$ was decomposed into squares of the same
mesh sizes as described in Section \ref{sec:hyb}, see Figure
\ref{fig:F1}.  The time step was chosen corresponding to the
stability criterion \eqref{CFL1}  as $\tau_l = 0.025
\cdot 2^{-l}$ for  $l=3,4,5,6$. Convergence results of the proposed finite
element scheme computed in $L_2$
and $H^1$ norms are presented in Tables \ref{testm4} - \ref{testm8}  for  $m=2,4,6,8$ in \eqref{eps}.
Relative norms in these tables were computed as
\begin{equation}\label{norms}
  \begin{split}
e_l^1 &= \displaystyle \frac{\displaystyle \max_{1 \leq k \leq N} \| E^k- E_h^k \|}{\displaystyle \max_{1 \leq k \leq N}\| E^k \|},\\
e_l^2 &= \displaystyle \frac{ \displaystyle \max_{1 \leq k \leq N} \| \nabla (E^k - E^k_h) \|}
{\displaystyle \max_{1 \leq k \leq N} \| \nabla {E^k} \|}.
\end{split}
  \end{equation}
$E$ and $E_h$ are  the exact and computed FE solutions  in  $\Omega_{FEM} \times (0,T)$, respectively, and  $N= T /\tau_l$.
 Logarithmic convergence  rates $r_1, r_2$   in these tables  are computed, viz. 
\begin{equation}
  \begin{split}
  r_1 &= \frac{ \left| \log \left( \frac{  el^1_h }{  el^1_{2h}} \right) \right|}{|\log(0.5)|},\quad\mbox{and}\quad 
  r_2 = \frac{\left | \log \left( \frac{  el^2_h}{ el^2_{2h}} \right) \right | }{|\log(0.5)|},
  \end{split}
  \end{equation}
where $el^{1,2}_h, el^{1,2}_{2h} $ are relative norms 
computed via \eqref{norms} on the mesh ${\mathcal T}_h$ with the mesh
size $h$ and $2h$, respectively.  Figure \ref{fig:F3} shows
convergence of the relative $L_2$ and $H^1$ norms computed via
\eqref{norms} and compared with exact behavior of $h$ and $h^2$.

\begin{table}[h!] 
\center
\begin{tabular}{ l  l l l l l l l l }
\hline
$l$ &  $nel$  & $nno$ &  $e_l^1$  &   $ \frac{  el^1_h }{  el^1_{2h}}$  &  $r_1$ & $e_l^2$  & $ \frac{  el^2_h }{  el^2_{2h}}$ & $r_2$ \\
\hline 
$3$ & $128$ & $81$    &  $4.878\cdot 10^{-2}$      & $ -$   &  $ -$ &  $3.902 \cdot 10^{-1} $  &  $-$      &   $ -$ \\
$4$ & $512$ & $289$   &  $1.222\cdot 10^{-2}$      &$3.992$ & $1.997$ &  $1.955 \cdot 10^{-1}$  &  $1.996$  & $0.997$ \\
$5$ & $2048$ & $1089$ &  $2.654\cdot 10^{-3}$       &$4.604$ & $2.203$ &  $8.492 \cdot 10^{-2}$   &  $2.302$ & $1.203$\\
$6$ & $8192$ & $4225$ &  $5.15\cdot 10^{-4}$        &$5.151$ & $2.365$ &  $3.297 \cdot 10^{-2}$  &  $2.575$     &  $1.365$ \\
\hline
\end{tabular}
\caption{Relative errors  $e_l^1$  and   $e_l^2$  in the $L_2$- and  $H^1$-norms, respectively, for
  mesh sizes $h_l= 2^{-l}, l=3,...,6$,  for $m=2$ in \eqref{eps}.}
\label{testm2}
\end{table}


\begin{table}[h!] 
\center
\begin{tabular}{ l  l l l l l l l l }
\hline
$l$ &  $nel$  & $nno$ &  $e_l^1$  &   $ \frac{  el^1_h }{  el^1_{2h}}$  &  $r_1$ & $e_l^2$  & $ \frac{  el^2_h }{  el^2_{2h}}$ & $r_2$ \\
\hline 
$3$ & $128$ & $81$    &  $5.54 \cdot 10^{-2}$      & $ -$      & $- $      &  $4.432 \cdot 10^{-1} $  &  $-$      &  $ -$ \\
$4$ & $512$ & $289$   &  $8.438 \cdot 10^{-3}$       & $6.565$ & $2.7148$ &  $ 1.35 \cdot 10^{-1}$  &  $3.283$  & $1.715$ \\
$5$ & $2048$ & $1089$ &  $1.581 \cdot 10^{-3}$       & $5.337$ & $2.4160$ & $ 5.06  \cdot 10^{-2}$   &  $2.668$ & $1.416$\\
$6$ & $8192$ & $4225$ &  $3.35 \cdot 10^{-4}$        & $4.722$ & $2.2394$ & $ 2.143 \cdot 10^{-2}$  &  $2.361$     &  $1.239$ \\
\hline
\end{tabular}
\caption{Relative errors  $e_l^1$  and   $e_l^2$  in the $L_2$-norm and in the $H^1$-norm, respectively, for
  mesh sizes $h_l= 2^{-l}, l=3,...,6$,  for $m=4$ in \eqref{eps}.}
\label{testm4}
\end{table}

\begin{table}[h!] 
\center
\begin{tabular}{ l  l l l l l l l l }
\hline
$l$ &  $nel$  & $nno$ &  $e_l^1$  &   $ \frac{  el^1_h }{  el^1_{2h}}$  &  $r_1$ & $e_l^2$  & $ \frac{  el^2_h }{  el^2_{2h}}$ & $r_2$ \\
\hline 
$3$ & $128$  & $81$    &  $ 1.856 \cdot 10^{-2}$  & $ -$     &$ -$       &  $  1.485 \cdot 10^{-1} $  &  $-$      &  $ -$ \\
$4$ & $512$  & $289$   &  $ 5.168 \cdot 10^{-3}$   & $ 3.592$ & $1.845$ &  $   8.268 \cdot 10^{-2}$  &  $1.796$  & $0.778$ \\
$5$ & $2048$ & $1089$  &  $1.594 \cdot 10^{-3}$    & $ 3.243$ & $1.697$ &  $  5.099 \cdot 10^{-2}$   &  $1.621$  & $0.697$\\
$6$ & $8192$ & $4225$  &  $ 3.6  \cdot 10^{-4}$    & $ 4.432$ & $2.148$ &  $ 2.301 \cdot 10^{-2}$    &  $2.216$  &  $1.148$ \\
\hline
\end{tabular}
\caption{Relative errors  $e_l^1$  and   $e_l^2$ in the $L_2$-norm and in the $H^1$-norm, respectively, for
  mesh sizes $h_l= 2^{-l}, l=3,...,6$,  for $m=6$ in \eqref{eps}.}
\label{testm6}
\end{table}

\begin{table}[h!] 
\center
\begin{tabular}{ l  l l l l l l l l }
\hline
$l$ &  $nel$  & $nno$ &  $e_l^1$  &   $ \frac{  el^1_h }{  el^1_{2h}}$  &  $r_1$ & $e_l^2$  & $ \frac{  el^2_h }{  el^2_{2h}}$ & $r_2$ \\
\hline 
$3$ & $128$  & $81$    &  $ 1.131 \cdot 10^{-2}$  & $ -$     &$ -$       &  $    1.045 \cdot 10^{-1} $  &  $-$      &  $ -$ \\
$4$ & $512$  & $289$   &  $ 5.669 \cdot 10^{-3}$   & $  2.304$ & $1.204$ &  $   9.071  \cdot 10^{-2}$  &  $1.152$  & $0.2$ \\
$5$ & $2048$ & $1089$  &  $ 1.711 \cdot 10^{-3}$    & $  3.314$ & $1.728$ &  $  5.475 \cdot 10^{-2}$   &  $1.657$  & $0.728$\\
$6$ & $8192$ & $4225$  &  $ 3.83  \cdot 10^{-4}$    & $ 4.468$ & $2.16$ &  $  2.451 \cdot 10^{-2}$    &  $2.234$  &  $1.16$ \\
\hline
\end{tabular}
\caption{Relative errors   $e_l^1$  and   $e_l^2$   in the $L_2$-norm and in the $H^1$-norm, respectively,
  for
  mesh sizes $h_l= 2^{-l}, l=3,...,6$,  for $m=8$ in \eqref{eps}.}
\label{testm8}
\end{table}
\begin{figure}[h!]
\begin{center}
  \begin{tabular}{cccc}
     \hline
      \multicolumn{4}{c}{ Computational meshes in the domain decomposition of     $\Omega = \Omega_{\rm FEM}\cup \Omega_{\rm FDM}$}\\
  {\includegraphics[scale=0.15, trim = 9.0cm 0.0cm 7.0cm 0.0cm, clip=]{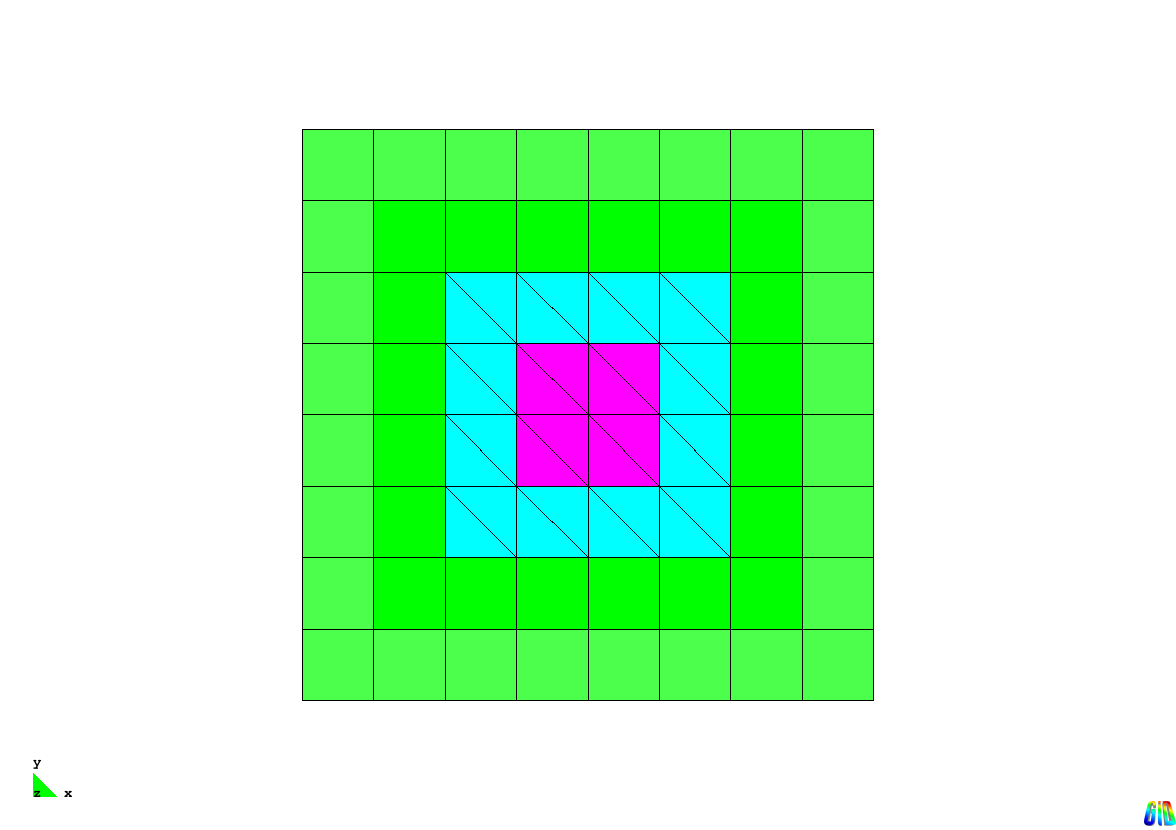}}  &
  {\includegraphics[scale=0.15, trim = 9.0cm 0.0cm 7.0cm 0.0cm, clip=]{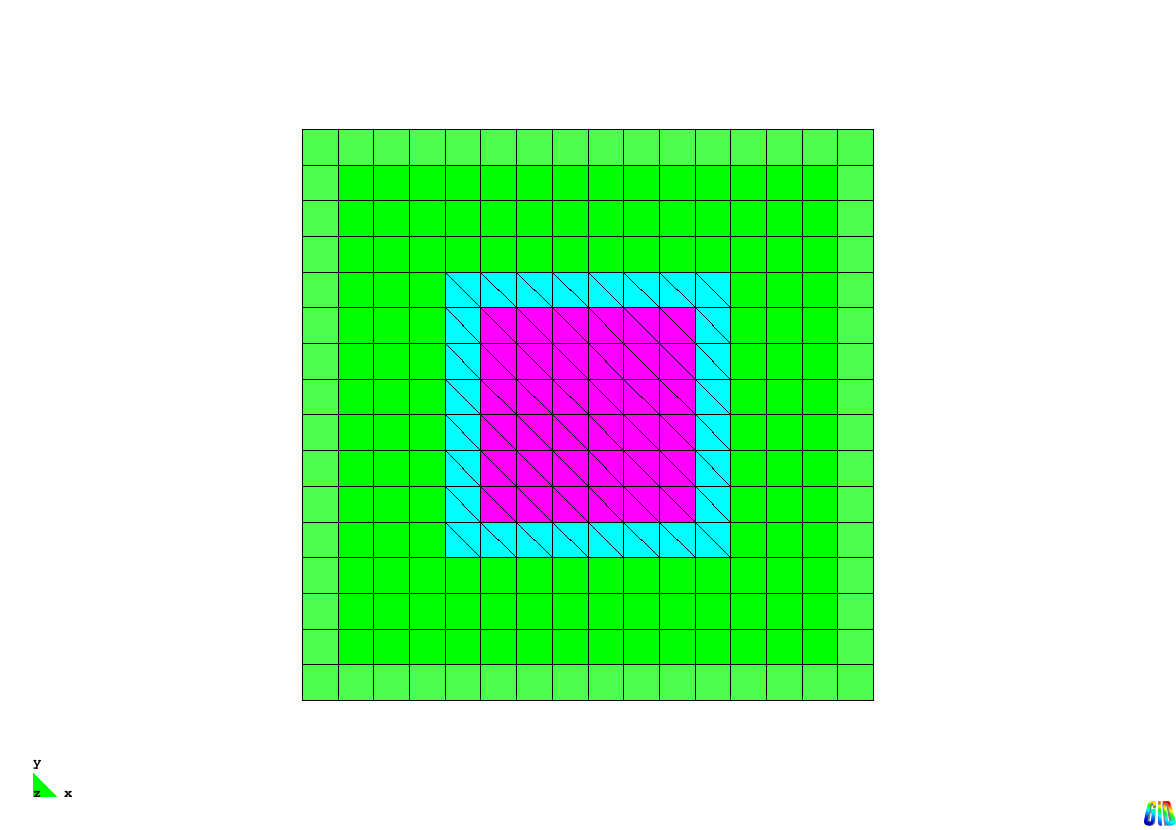}} &
  {\includegraphics[scale=0.15,trim = 9.0cm 0.0cm 7.0cm 0.0cm, clip=]{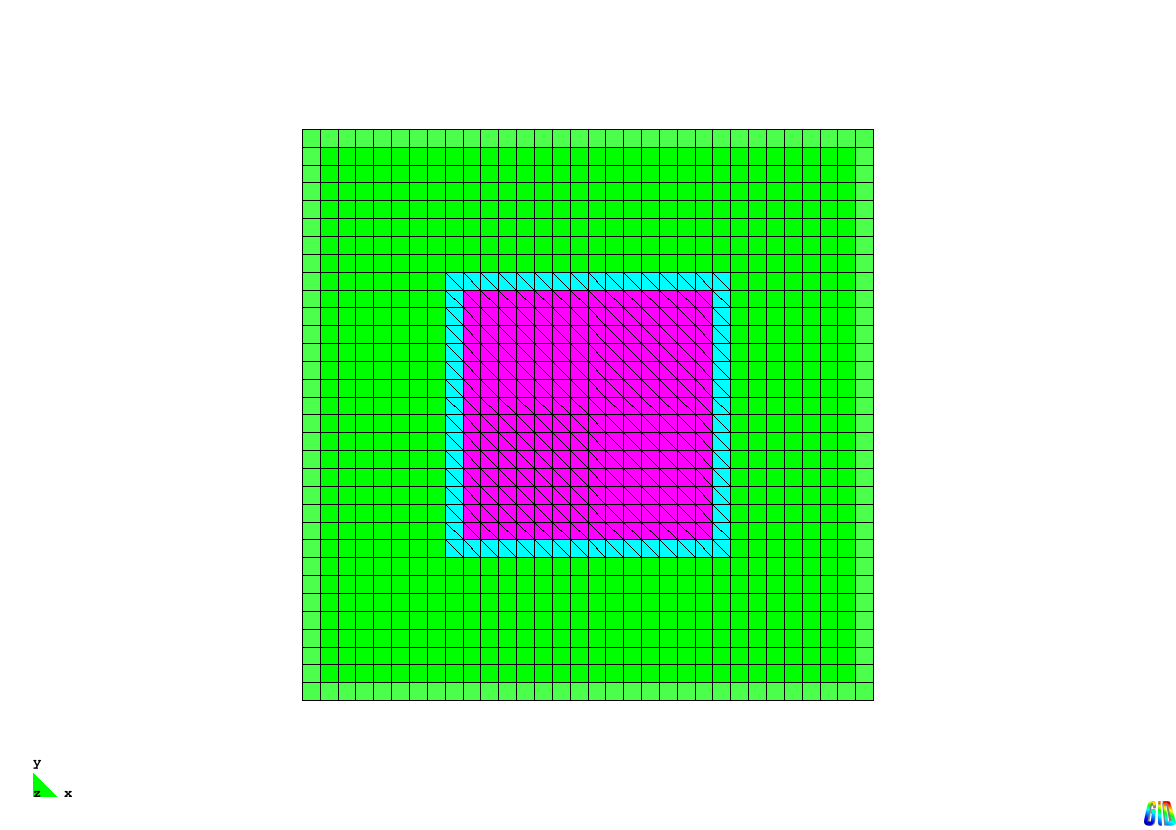}}  &
   {\includegraphics[scale=0.15,trim = 9.0cm 0.0cm 7.0cm 0.0cm, clip=]{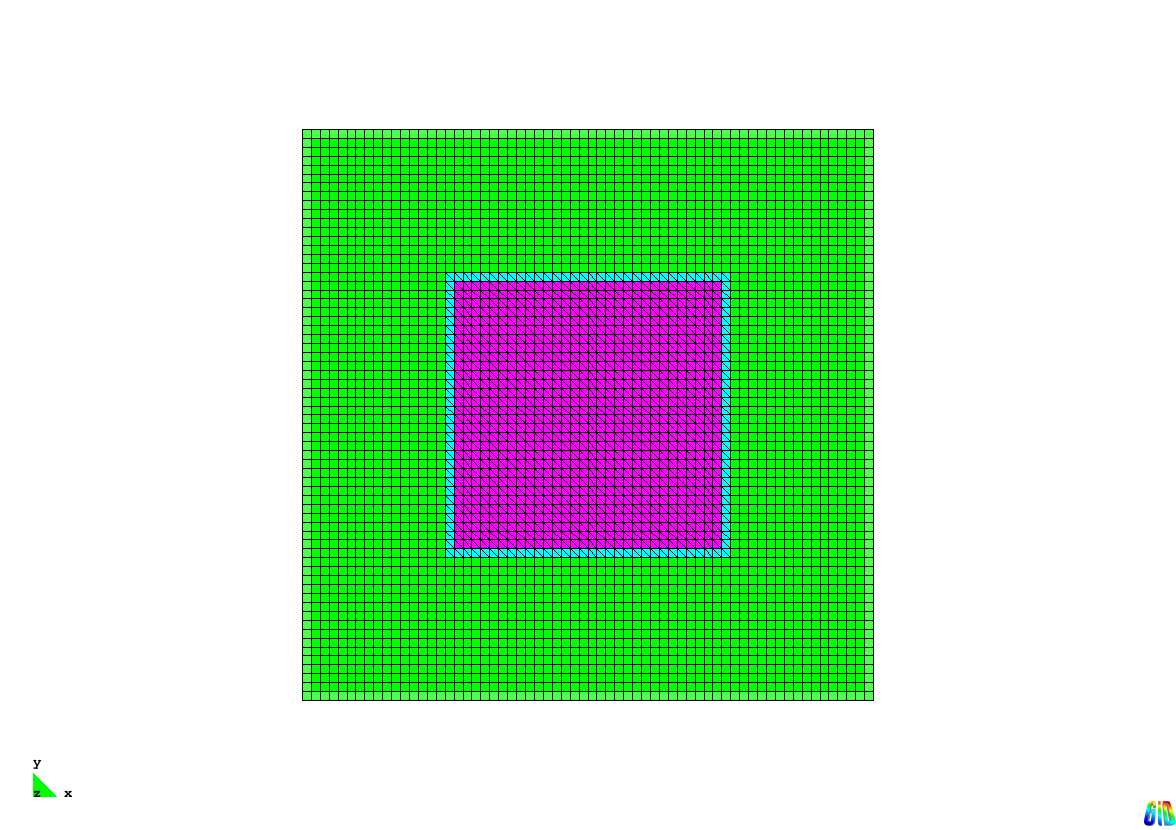}}\\
    $h=0.125$ &  $h=0.0625$ & $h=0.03125$ & $h=0.015625$\\
    \hline
      \multicolumn{4}{c}{ Exact solution   $|E|, m=8$, in $\Omega$} \\
  {\includegraphics[scale=0.1, trim = 6.0cm 0.0cm 1.0cm 0.0cm, clip=]{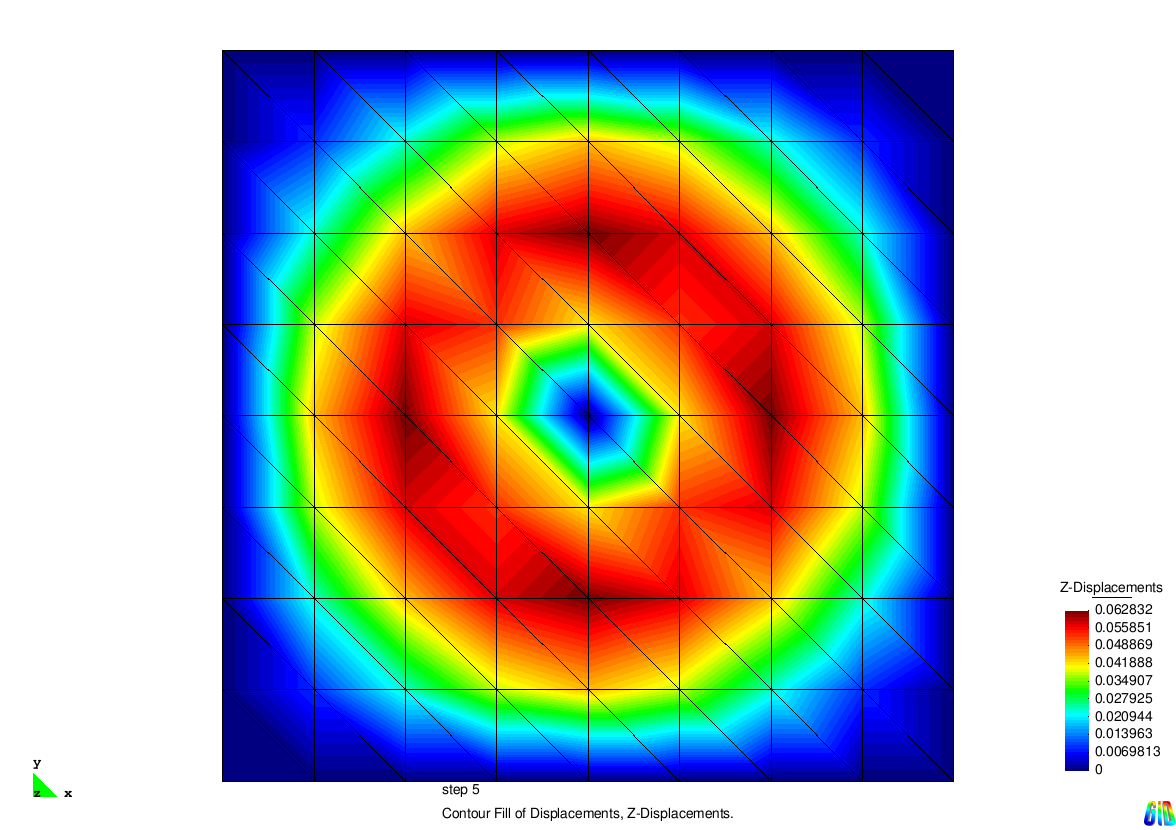}}  &
  {\includegraphics[scale=0.1, trim = 6.0cm 0.0cm 1.0cm 0.0cm, clip=]{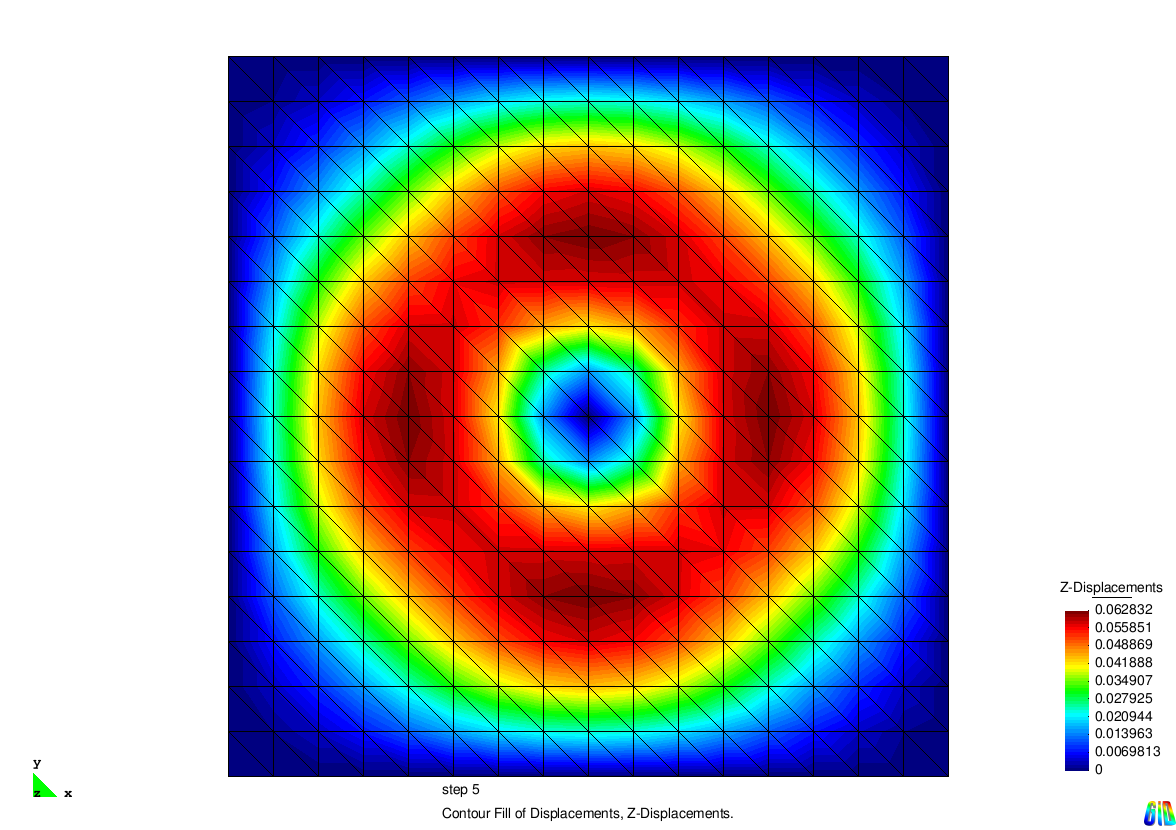}} &
  {\includegraphics[scale=0.1,trim = 6.0cm 0.0cm 1.0cm 0.0cm, clip=]{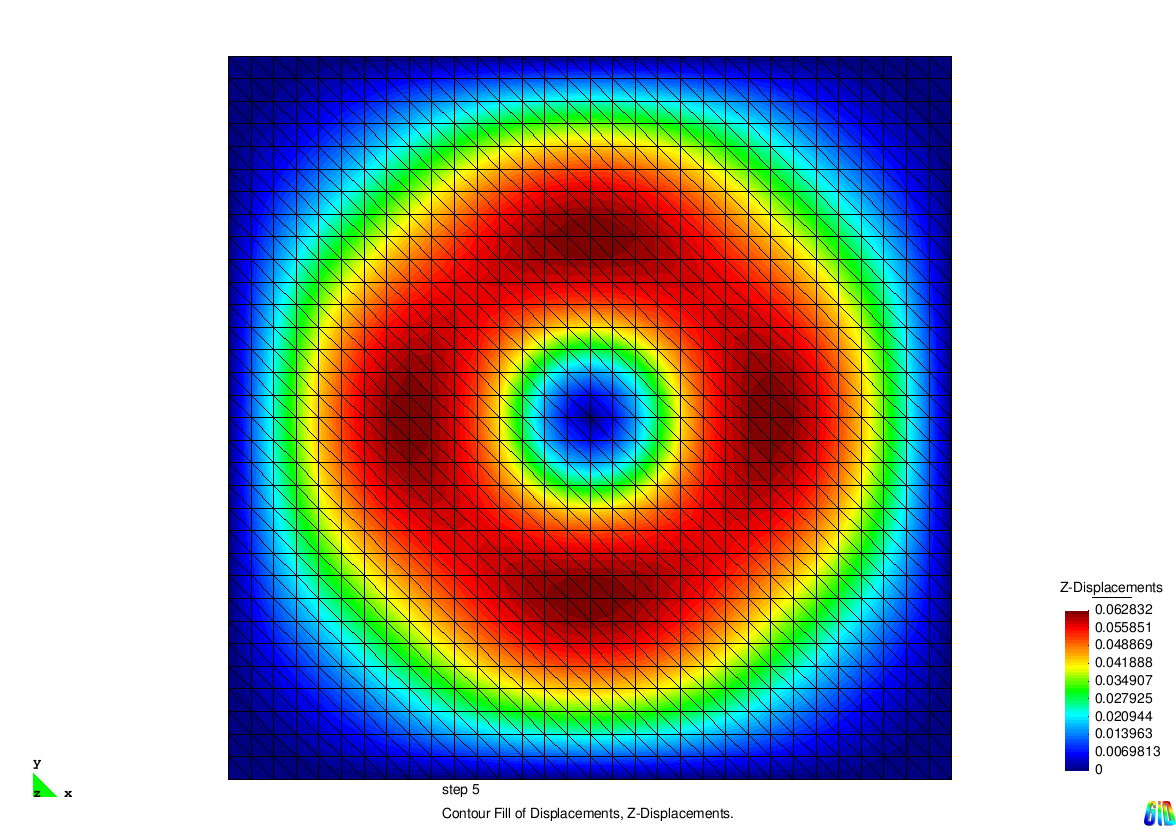}}  &
  {\includegraphics[scale=0.1,trim = 6.0cm 0.0cm 1.0cm 0.0cm, clip=]{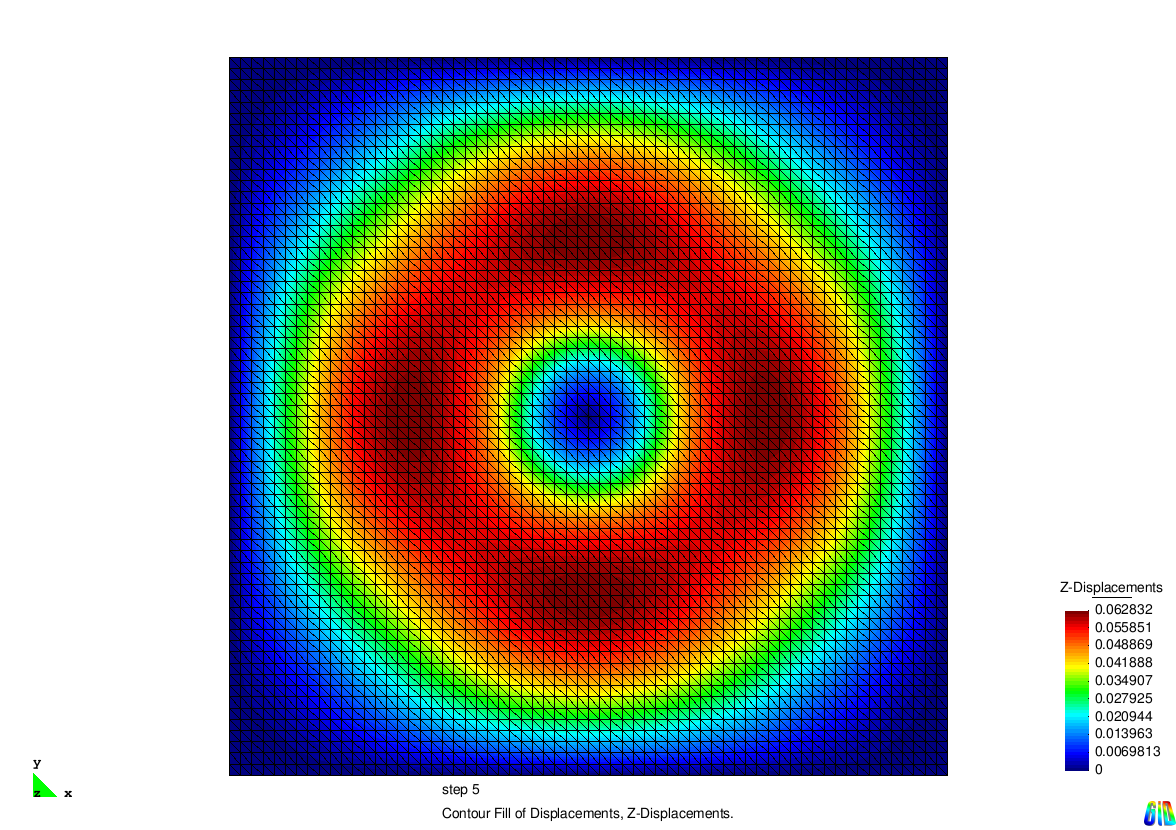}}\\
    $h=0.125$ &  $h=0.0625$ & $h=0.03125$ & $h=0.015625$\\
   \hline
   \multicolumn{4}{c}{Computed domain decomposition solution   $|E_h|, m=8$,  in $\Omega = \Omega_{\rm FEM}\cup \Omega_{\rm FDM}$ }\\
     {\includegraphics[scale=0.1, trim = 6.0cm 0.0cm 1.0cm 0.0cm, clip=]{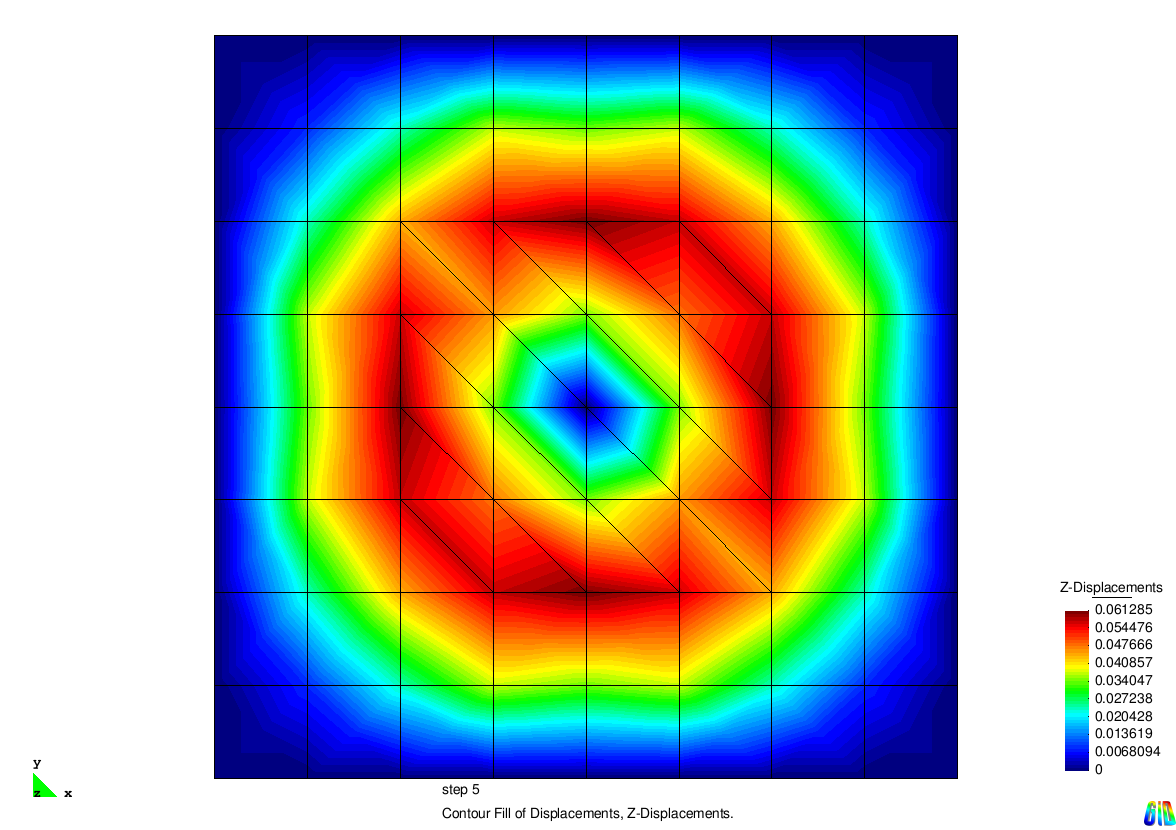}}  &
  {\includegraphics[scale=0.1,  trim = 6.0cm 0.0cm 1.0cm 0.0cm, clip=]{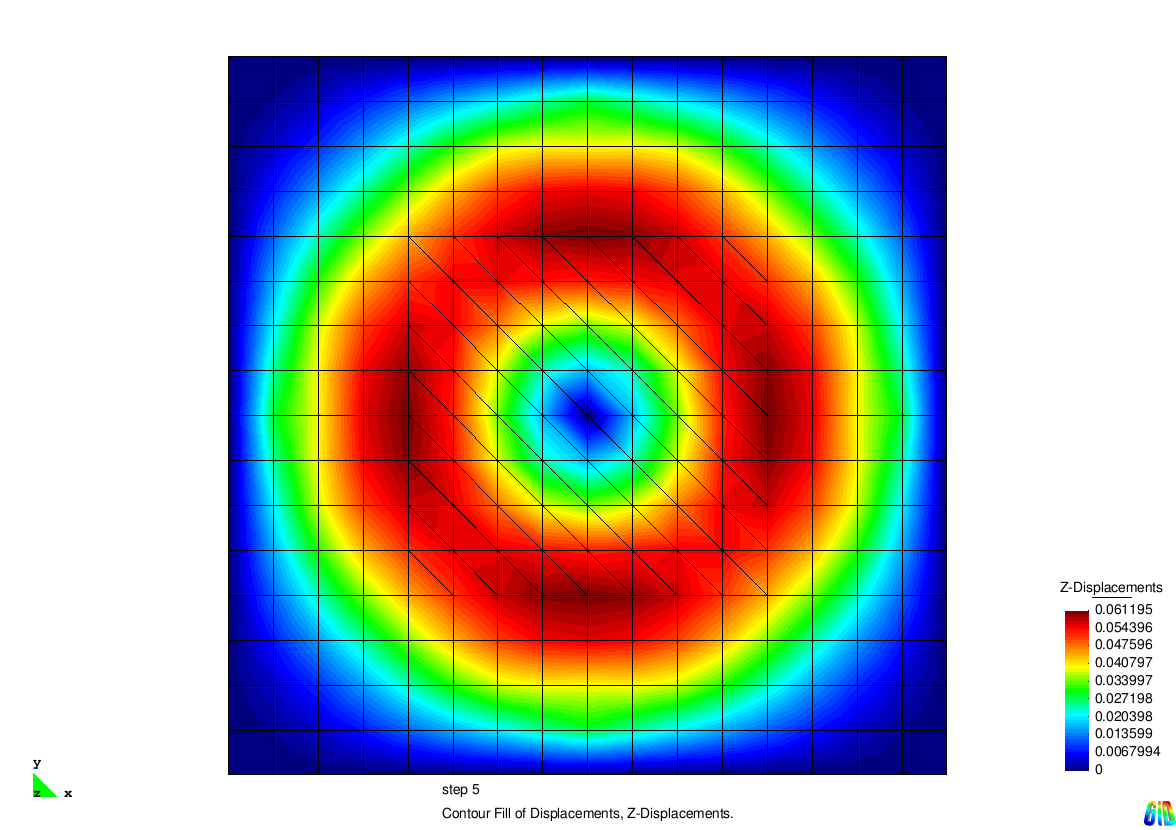}} &
  {\includegraphics[scale=0.1, trim = 6.0cm 0.0cm 1.0cm 0.0cm,  clip=]{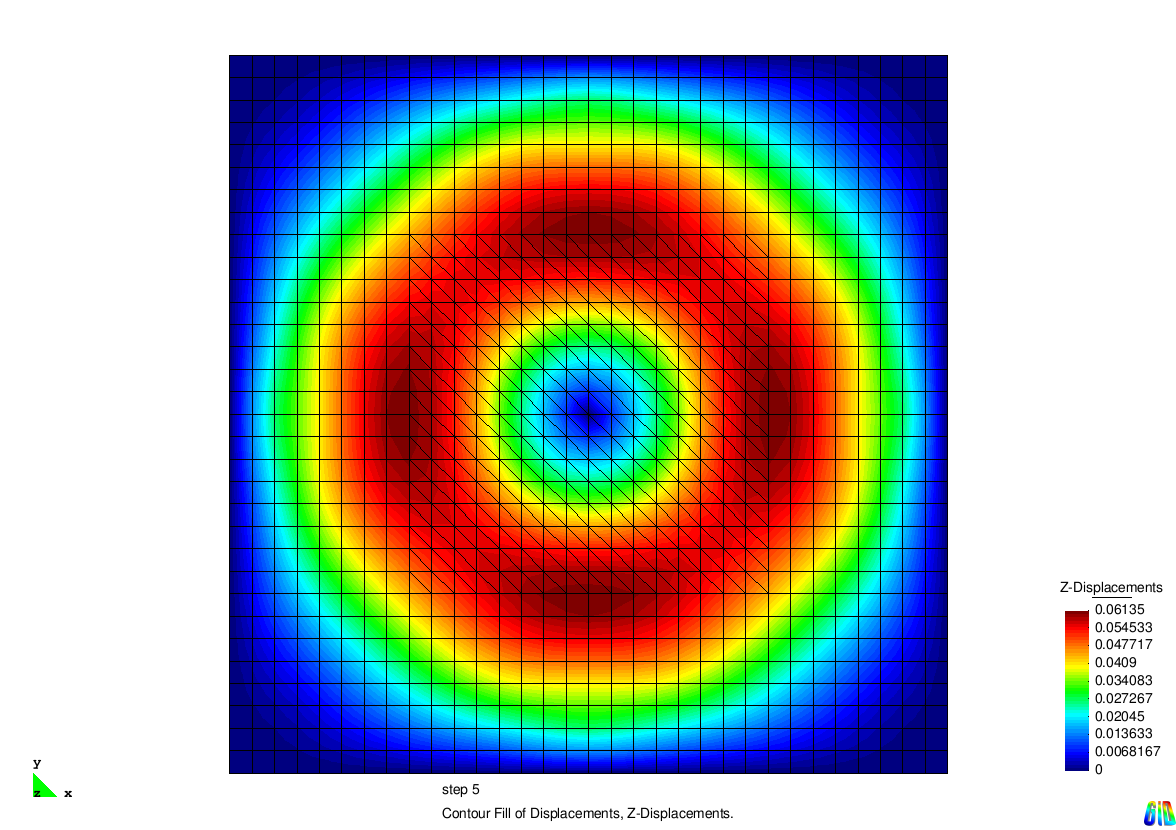}}  &
   {\includegraphics[scale=0.1, trim = 6.0cm 0.0cm 1.0cm 0.0cm,  clip=]{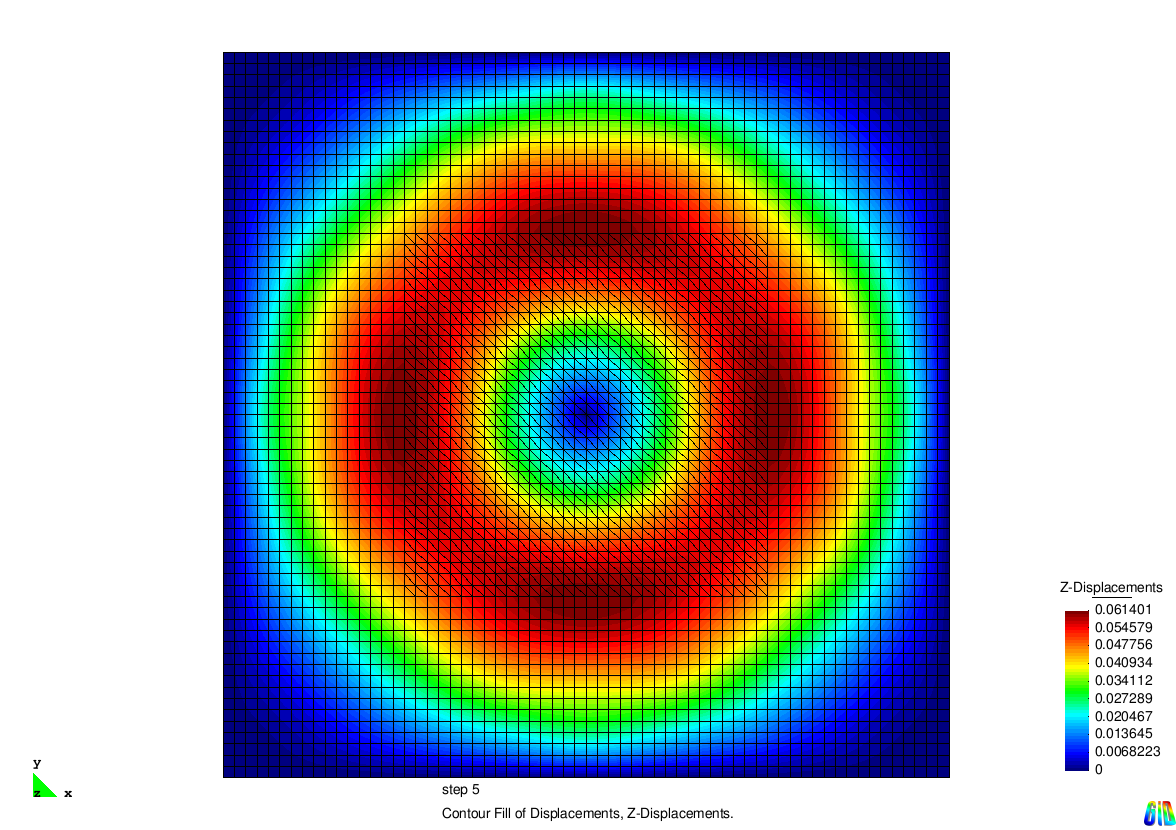}}\\
     $h=0.125$ &  $h=0.0625$ & $h=0.03125$ & $h=0.015625$\\ 
   \hline
       \multicolumn{4}{c}{Computed finite difference solution  $|\hat{E}_h|, m=8$,  in the domain decomposition algorithm   in $\Omega_{\rm FDM}$}\\
     {\includegraphics[scale=0.1, trim = 6.0cm 0.0cm 1.0cm 0.0cm, clip=]{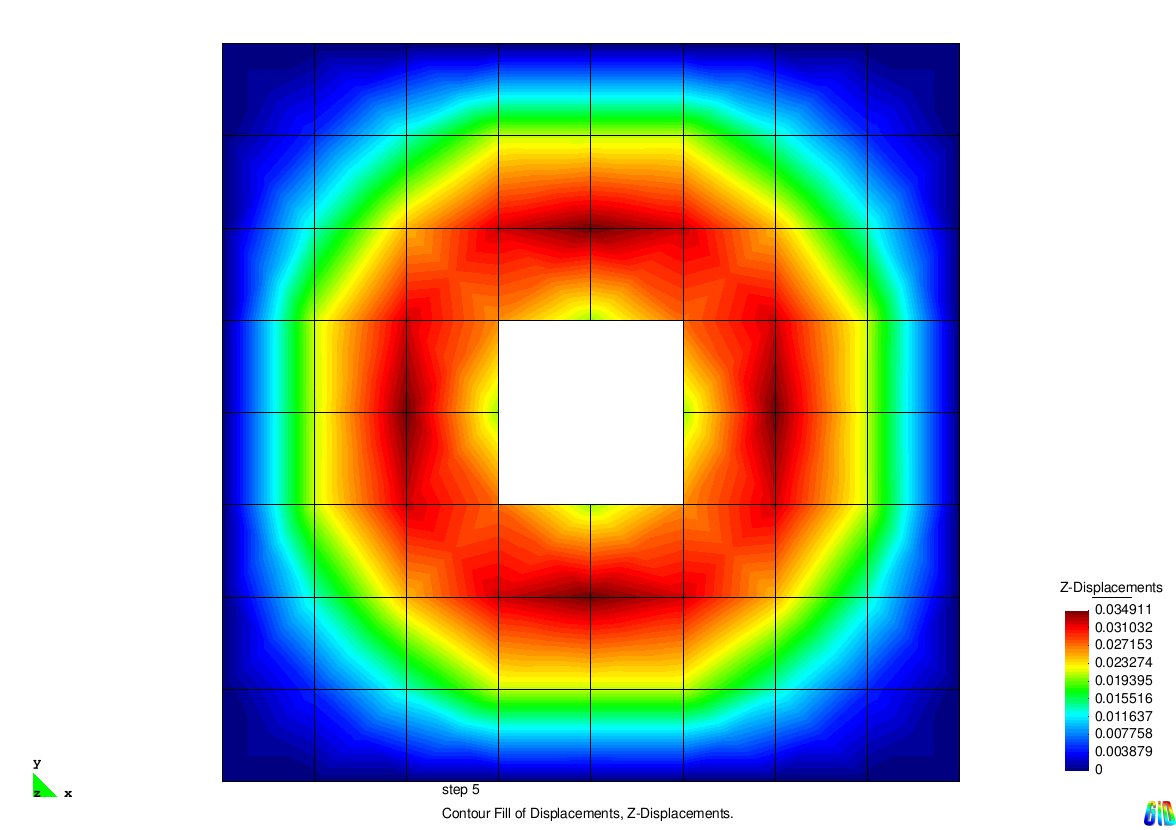}}  &
  {\includegraphics[scale=0.1,  trim = 6.0cm 0.0cm 1.0cm 0.0cm, clip=]{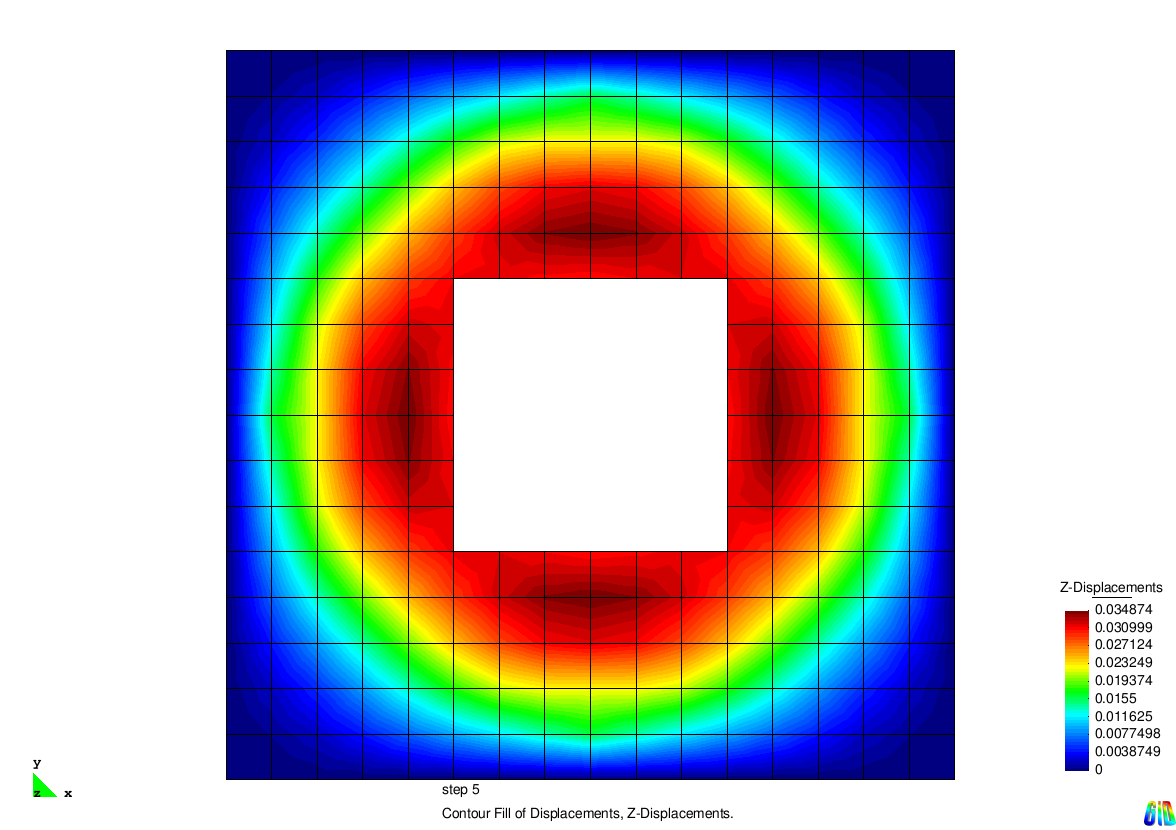}} &
  {\includegraphics[scale=0.1, trim = 6.0cm 0.0cm 1.0cm 0.0cm,  clip=]{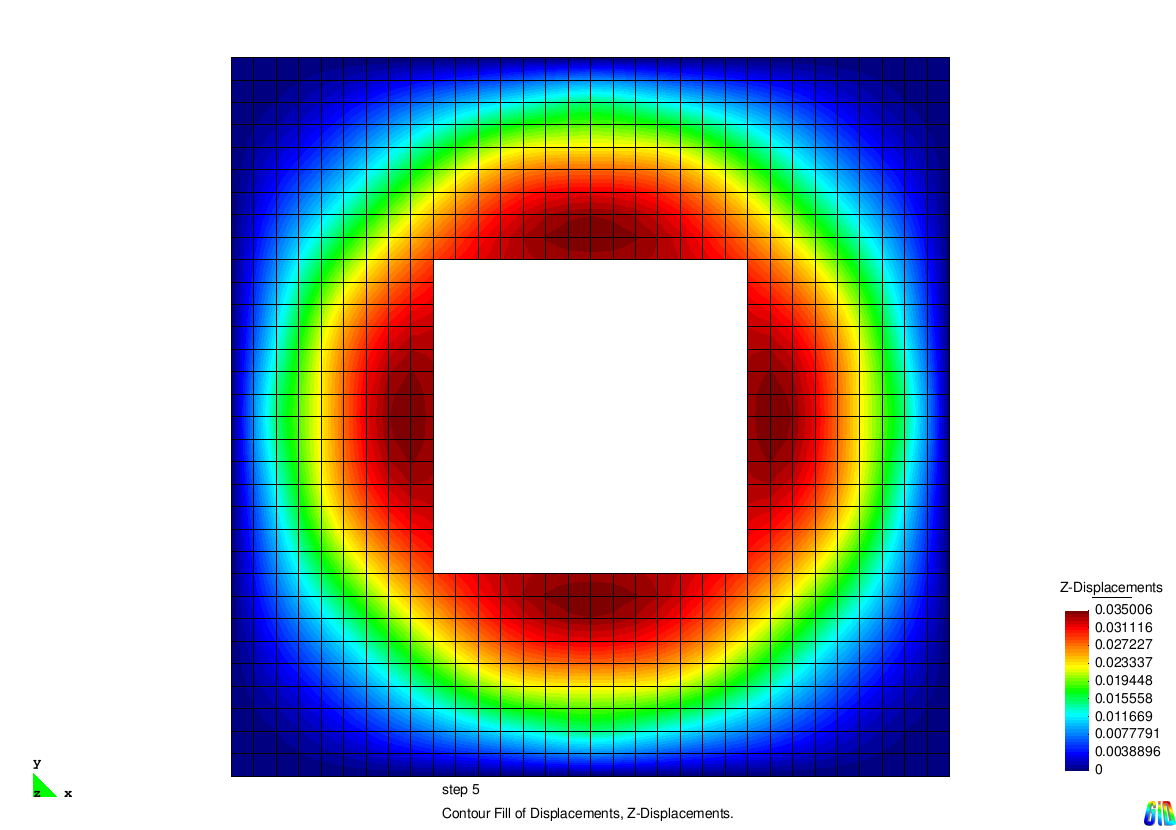}}  &
   {\includegraphics[scale=0.1, trim = 6.0cm 0.0cm 1.0cm 0.0cm,  clip=]{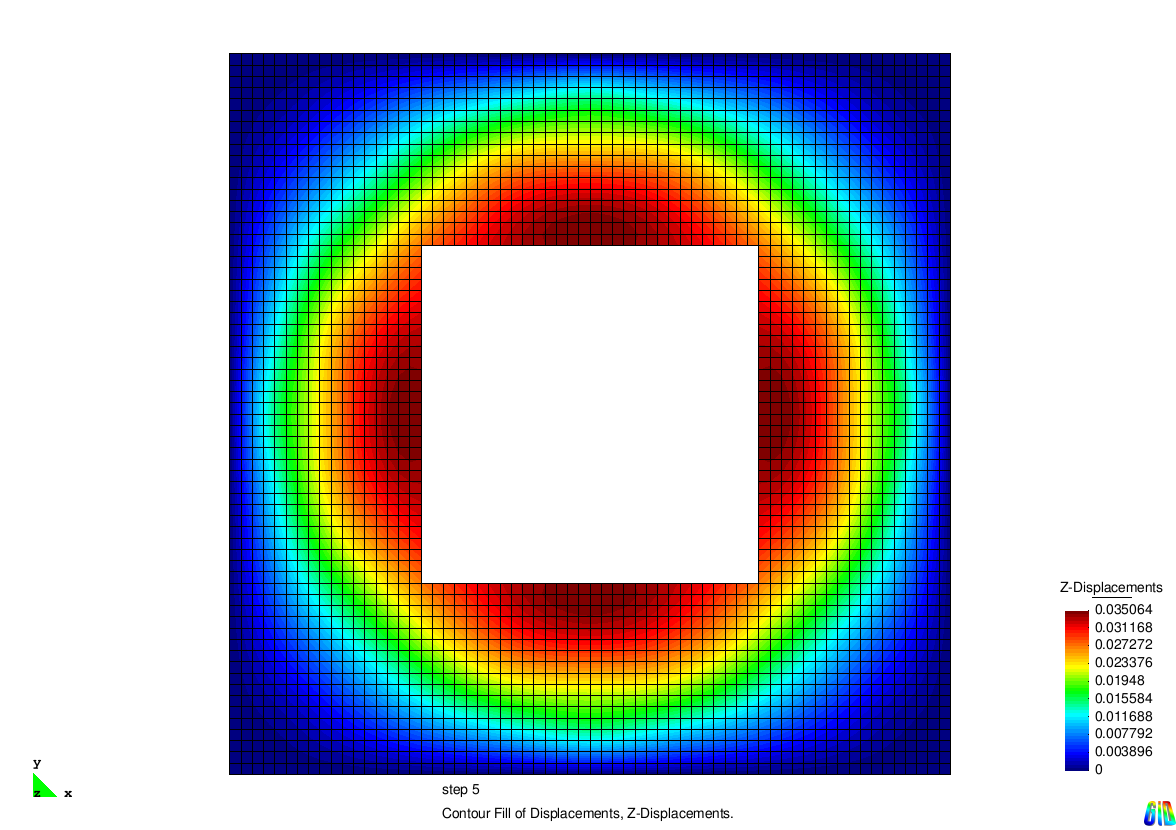}}\\
     $h=0.125$ &  $h=0.0625$ & $h=0.03125$ & $h=0.015625$\\ 
   \hline
\end{tabular}
\end{center}
\caption{
    \emph{Computed vs. exact solution at the time $t=0.25$ for different meshes taking $m=8$    in \eqref{eps}. Algorithm 2 was used in the domain decomposition method. Common  elements in  $\Omega_{\rm FEM}$  and $\Omega_{\rm FDM}$ on different meshes are presented on the top figures  and are outlined by light blue color.  We observe smooth  hybrid solution across FE/FD boundaries. }}
\label{fig:F6}
\end{figure}

Computed hybrid FE/FD versus exact solutions with $m=8$ in \eqref{eps}
at the time $t=0.25$ are presented in Figures \ref{fig:F6}, where 
$|E_h|$ is computed 
using the domain decomposition algorithm (Algorithm 2)
for different meshes with sizes $h_l= 2^{-l}, l= 3,4,5,6$. The top figures
 of Figure \ref{fig:F6}
 present hybrid FE/FD meshes which were used for computations;
 common hybrid FE/FD solution in $\Omega$ is presented in the middle figures, and bottom figures
show only FD solution as the part of the common hybrid solution in
$\Omega_{\rm FDM}$.    Interpreting these
figures we observe smooth behavior of the hybrid solution across
finite element/finite difference boundary, as was predicted in theory.

Furthermore, through these figures, as well as tables and Figure \ref{fig:F3}, we
observe that with increasing $l$ in $h_l= 2^{-l}, l=3,4,5,6$, the 
computational errors approach the second order convergence in
$L_2$- and first order in $H^1$-norm for
$m=2,4,6,8$.  Therefore, we can conclude that the finite element
scheme in the hybrid FE/FD method; considered in $\Omega_{\rm FEM} $, behaves like a first order
method in $H^1$-norm and a second order method in the  $L_2$-norm.
These results are all in good agreement with the analytic estimates
derived in Sections 5-6, as well as with results presented for 
finite element method in \cite{BR1, BR2} for the whole $\Omega$.

\vskip 0.5cm 

\section*{Conclusion}

In this paper we  present   stability and convergence analysis
for the domain decomposition
FE/FD method for time-dependent 
Maxwell's equations  developed in \cite{BGrote, BMaxwell}.
The convergence is optimal due to the assumed maximal available regularity of the exact solution 
in a Sobolev space. 

The analysis are performed for the semi-discrete (spatial discretization) problem for, the   constructed, finite 
element schemes in two different 
settings: in $\Omega$ and $\Omega_{FEM}$. The temporal discretization algorithms are constructed using 
the CFL condition \eqref{CFL1} derived in  \cite{BR1}. 
We have  
implemented several numerical 
examples that validate the robustness of the theoretical studies. 

In a forthcoming, complementary, study we plan to extend the results in here to a problem with the 
presence of electrical conductivity term $\sigma\partial_t E$ which renders the equation to an parabolic-hyperbolic one.

\section*{Acknowledgments}

The research of both authors
is supported by the Swedish Research Council grant VR 2018-03661. The first 
author acknowledges the support of the VR grant DREAM.

\end{document}